\documentclass[3p, american, preprint, longtitle]{elsarticle}
\pdfoutput=1
\usepackage[utf8]{inputenc}
\usepackage[T1]{fontenc}
\usepackage{babel, lmodern}
\usepackage[babel=true]{microtype}
\usepackage{amsmath, amssymb, physics}
\usepackage{enumitem, subcaption, xcolor}
\usepackage{float}

\usepackage{hyperref}
\usepackage[capitalize]{cleveref} 
\crefname{equation}{}{}
\crefname{figure}{Fig.}{Figs.}
\Crefname{figure}{Figure}{Figures}
\crefname{section}{Sect.}{Sects.}
\Crefname{section}{Section}{Sections}

\biboptions{longnamesfirst}
\DeclareMathOperator{\conv}{conv}
\DeclareMathOperator{\vol}{vol}
\newcommand{\R}{\mathbb{R}}

\makeatletter
\def\ps@pprintTitle{%
   \def\@oddfoot{\reset@font\phantom{\today}\hfil\thepage\hfil\today}
}
\makeatother

\begin{document}

\begin{frontmatter}

\journal{Comput.\ Math.\ Math.\ Phys.}
\title{Non-simplicial Delaunay meshing via approximation by radical partitions}

\author[CCAS,MIPT]{Vladimir Garanzha}
\ead{garan@ccas.ru}
\ead[url]{http://www.ccas.ru/gridgen/lab}

\author[CCAS,MIPT]{Liudmila Kudryavtseva}
\ead{liukudr@yandex.ru}

\author[MIPT]{Lennard Kamenski}
\ead{l.kamenski@arcor.de}
\ead[url]{https://gitlab.com/lkamenski}

\address[CCAS]{Dorodnicyn Computing Center FRC CSC RAS, Moscow, Russia}
\address[MIPT]{Moscow Institute of Physics and Technology, Dolgoprudnyi, Russia}

\begin{abstract}
We consider the construction of a polyhedral Delaunay partition as a limit of the sequence of power diagrams (in Russian traditionally called \emph{radical partitions}), while the dual Voronoi diagram is obtained as a limit of the sequence of weighted Delaunay partitions.
Using a lifting analogy, this problem is reduced to the construction of a pair of dual convex polyhedra, inscribed and superscribed around a circular paraboloid, as a limit of the sequence of pairs of general dual convex polyhedra.
The sequence of primal polyhedra should converge to the superscribed polyhedron, while the sequence of dual polyhedra converges to the inscribed polyhedron.

Different rules can be used to build sequences of dual polyhedra.
We are mostly interested in the case when the vertices of primal polyhedra can move or merge together, meaning that no new faces are allowed for dual polyhedra.
These rules essentially define the transformation of the set of initial spheres defining a power diagram into the set of Delaunay spheres using sphere movement, radius variation, and sphere elimination as admissible operations.
Although a rigorous foundation (existence theorems) for this problem is still unavailable, we suggest a functional measuring the deviation of the convex polyhedron from the one inscribed into the paraboloid.
This functional is the discrete Dirichlet functional for the power function (power of a point with respect to a ball) which is a linear interpolant of the vertical distance of the dual vertices from the paraboloid.
The absolute minimizer of this functional is attained on the constant power field, meaning that the inscribed polyhedron can be obtained using a simple translation.
This formulation of the Dirichlet functional for the dual surface is not quadratic since the unknowns are the vertices of the primal polyhedron, hence, the transformation of the set of spheres into Delaunay spheres is not unique.

In this work, we concentrate on the experimental confirmation of the viability of the suggested approach and put aside mesh quality problems.
The zero value of the gradient of the proposed functional defines a manifold describing the evolution of Delaunay spheres.
Hence, Delaunay-Voronoi meshes can be optimized using this manifold as a constraint.
Numerical examples illustrate polygonal Delaunay meshing in planar domains.
\end{abstract}

\begin{keyword}
power diagram \sep{} radical partition
\sep{} Delaunay triangulation \sep{} weighted Delaunay triangulation
\sep{} Delaunay partition \sep{} Voronoi triangulation
\MSC[2020] 65M50 \sep{} 65N50 \sep{} 65L50
\\
\medskip{}
\footnotesize{}
This is a preprint of the Work accepted for publication in CMMP.
\copyright{} Pleiades Publishing, Ltd., 2022. \url{https://pleiades.online}
\end{keyword}

\end{frontmatter}

\section{Introduction}

In his famous talk ``\foreignlanguage{french}{Sur la sphere vide}'' at the 1924 Toronto Geometry Congress~\cite{Delone-1924} and later in his work~\cite{Delone-1936}, Boris Nikolaevich Delaunay (Delone) %
introduced a simplicial spatial partition of a given discrete vertex set consisting of $d$-dimensional simplexes with empty circumspheres, empty in the sense that they do not contain any vertices  inside.
The famous ``Delaunay Lemma'' states that if the empty circumsphere condition holds locally for any two adjacent simplexes having a common $d-1$-dimensional face, then all circumspheres are empty.
To prove this, Delaunay used the concept of the power $\tau$ of the point $\vb*{a}$ with respect to the ball $B$ of radius $R$ with center $\vb*{c}$,
\[
  \tau_B(\vb*{a}) = \qty|\vb*{c} - \vb*{a}|^2 - R^2
  .
\]

\begin{figure}[h]%
  \includegraphics[width=0.49\textwidth]{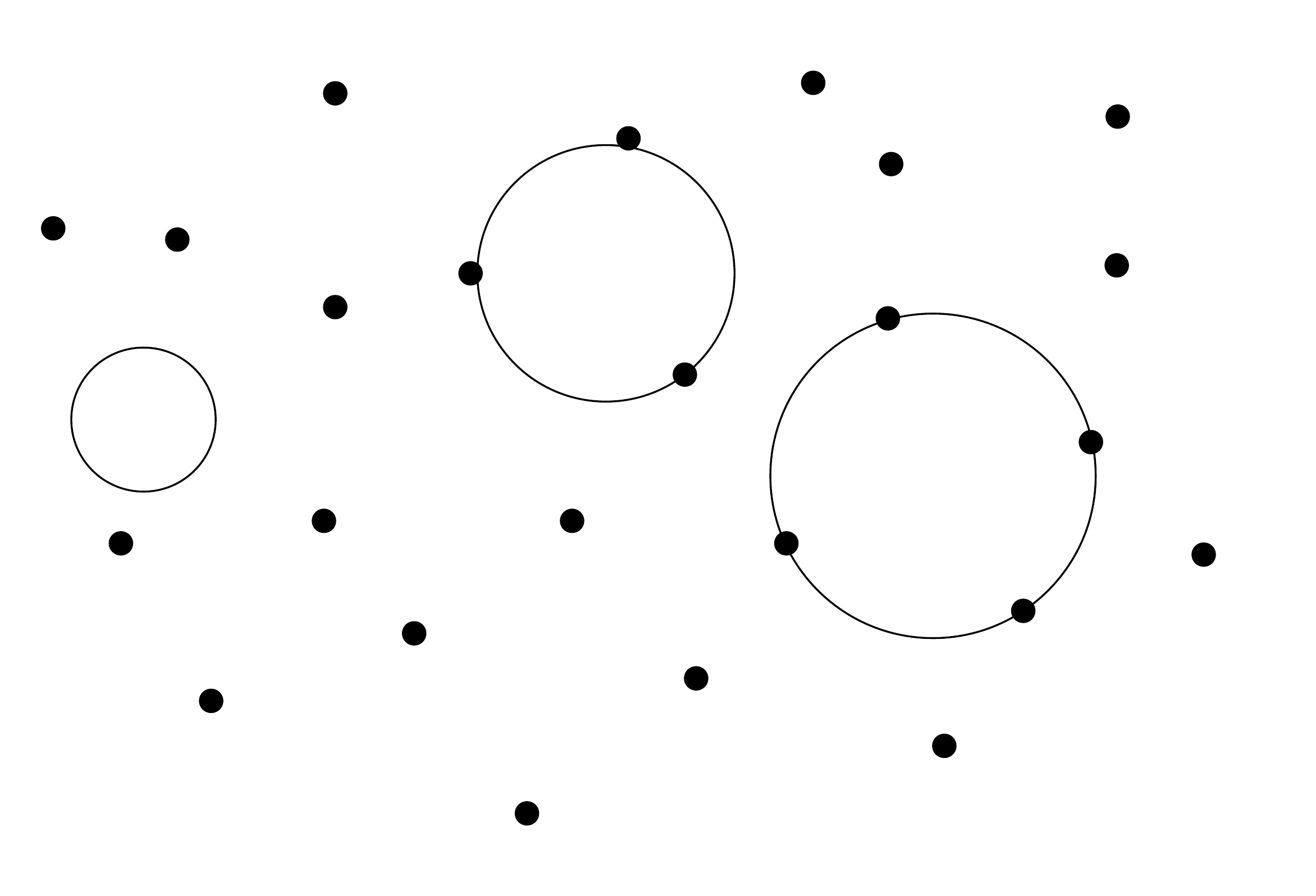}%
  \hfill{}%
  \includegraphics[width=0.49\textwidth]{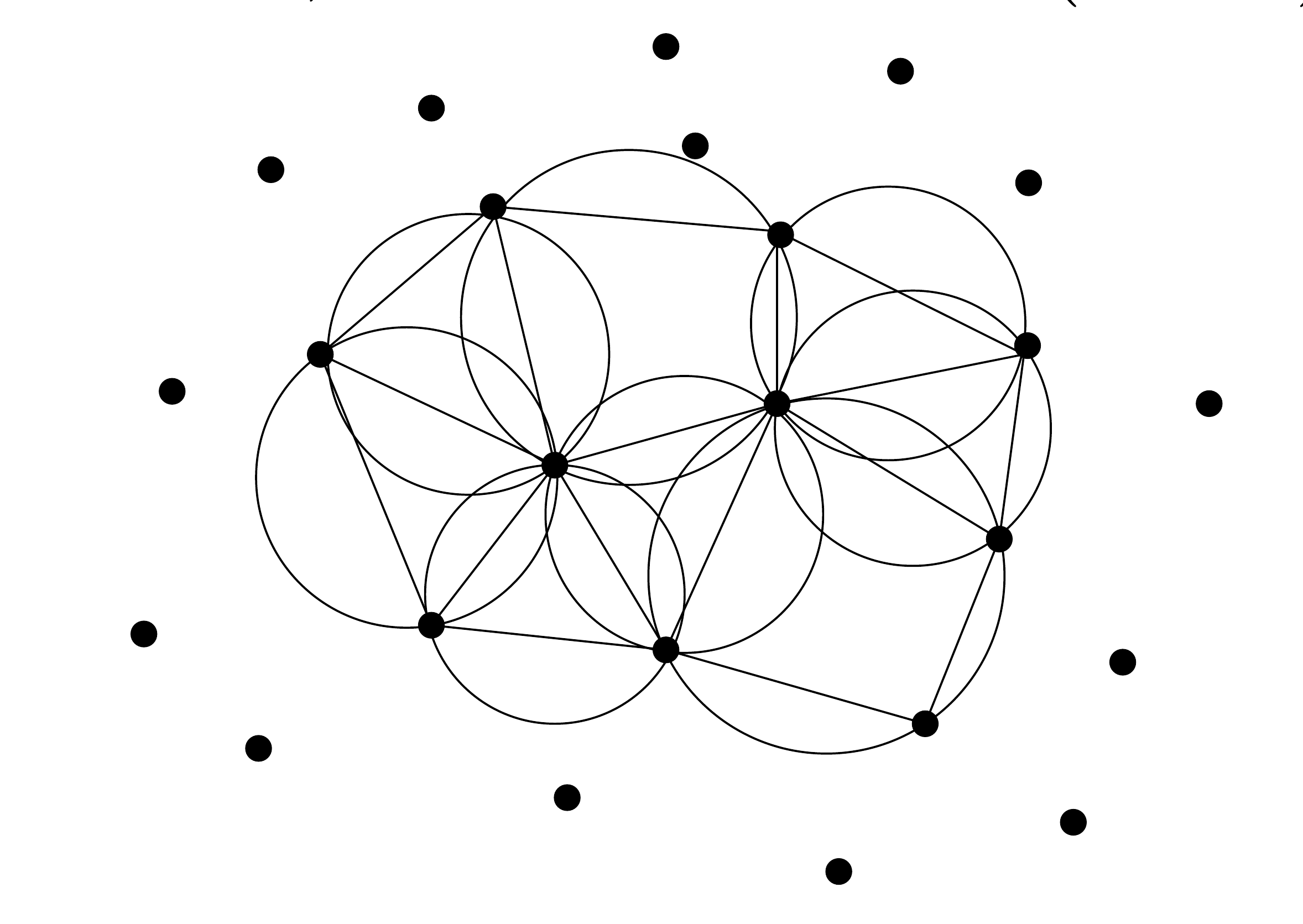}%
  \caption{%
     L-partition: (left) an empty sphere moving through the point set,
     (right) spheres with the $d$-dimensional point sets}%
  \label{fig:Delone:def}%
\end{figure}

In 1937, Delaunay introduced the \emph{L-partition}~\cite{Delone-1937}, which is obtained by moving, contracting, and inflating an empty sphere placed among the vertices of a discrete set in $\R^d$ (\cref{fig:Delone:def}, left).
In such a way, all possible empty spheres with a $d$-dimensional set of vertices on their surface, the so-called \emph{L-spheres}, are identified.

The convex envelope of the set of vertices lying on an L-sphere is a convex polyhedron called by Delaunay an \emph{L-polyhedron} (\cref{fig:Delone:def}, right).
The set of L-polyhedra defines a normal partition of the space called the L-partition.
The Delaunay lemma is naturally generalized for this polyhedral partition.
These days, the L-partition is called the \emph{Delaunay partition}, while the L-sphere, i.e., empty spheres with a $d$-dimensional set of vertices on their surfaces, are called the \emph{Delaunay spheres}.
The convex envelope of the centers of spheres passing through a vertex defines a Voronoi polyhedron for this vertex.
Note that both Delaunay and Voronoi partitions for a prescribed vertex set are unique.
The duality relations between Delaunay and Voronoi partitions are simple and elegant, which is not the case for Delaunay triangulations.
For example, for a set of co-circular vertices, the same circumcenter (Voronoi vertex) can be generated by several Delaunay simplexes, whereas the correspondence between a Delaunay polyhedron and the corresponding Voronoi vertex is always one-to-one.

Each Delaunay polyhedron can be split into simplexes, creating a simplicial partition generally called a \emph{Delaunay triangulation}, provided that the triangulations of the individual Delaunay polyhedra are consistent with each other.
Note the subtle but significant difference: a Delaunay edge/face in a Delaunay partition has at least one \emph{closed} empty circumball, whereas in a Delaunay triangulation it has at least one \emph{open} circumball.
Edges/faces of a triangulation which do not have \emph{closed} empty circumballs are the edges/faces added by triangulating Delaunay polyhedra.
These additional edges do not have a one-to-one correspondence with the dual Voronoi faces and the additional faces, correspondingly, do not have a one-to-one correspondence with the dual Voronoi edges.

Note that the construction based on a correct triangulation of a Delaunay partition excludes flat degenerate simplexes, the so-called \emph{slivers}.
Strictly speaking, a sliver is a degenerate simplex with a $d-m$-dimensional vertex set ($m>0$) for which it is possible to construct a circumscribed sphere of a finite radius.
Obviously, such a sphere is not unique.
Each such simplex is generated by an incorrect triangulation of a $d-m$-dimensional face of the Delaunay partition.
The sliver identification problem is aggravated if approximate computations are involved, which is, for example, the case of the standard floating point arithmetics.

\begin{figure}[t]%
   \includegraphics[width=0.49\textwidth]{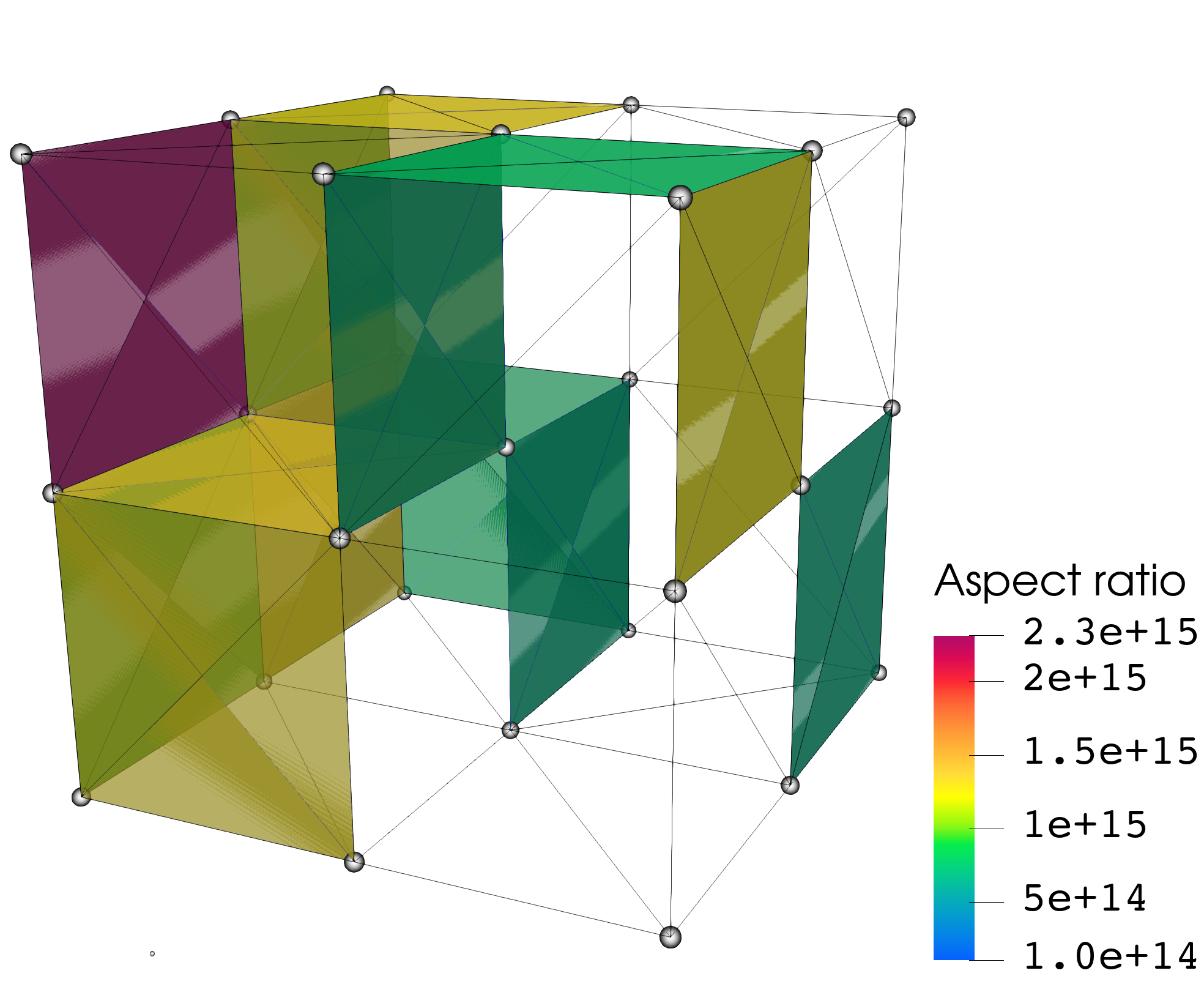}%
  \hfill{}%
  \includegraphics[width=0.49\textwidth]{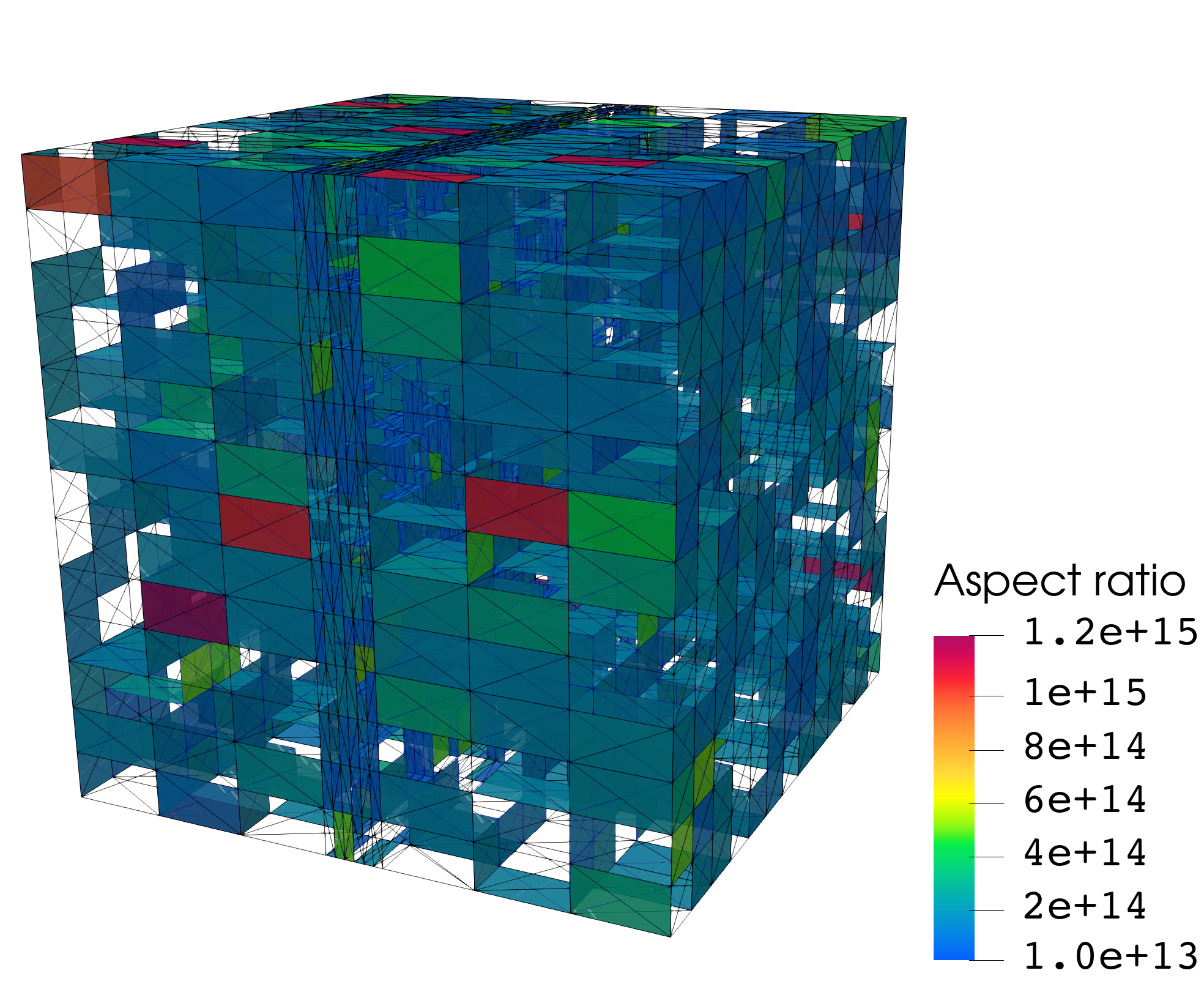}%
  \caption{Artificial slivers in \emph{TetGen} Delaunay meshes
    caused by rounding in the input data for
    a simple cubic lattice (left) and a more complex example (right)}%
  \label{fig:tetgen:slivers}
\end{figure}

\begin{figure}[t]%
  \includegraphics[width=0.48\textwidth]{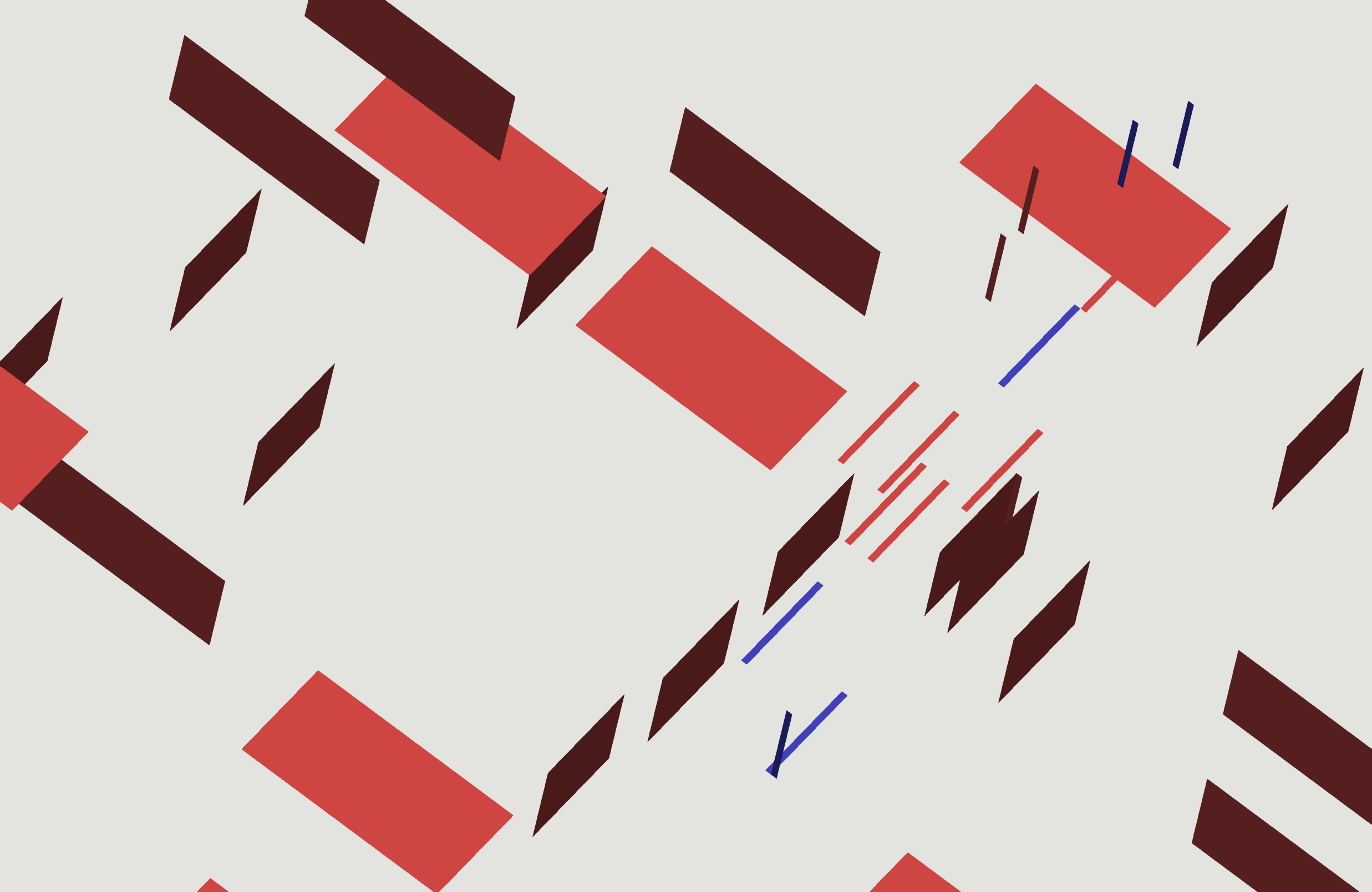}%
  \hfill{}%
  \includegraphics[height=0.48\textwidth,angle=90]{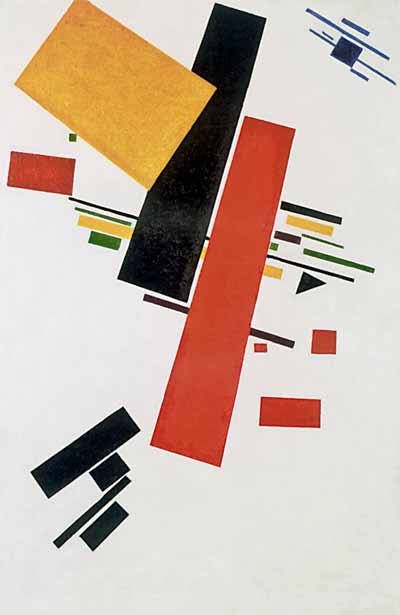}%
  \caption{``Sliver Art \textnumero~3'':
     with a proper chosen color map,
     the visualization of a selection of slivers (left) resembles
     the Supremus series of Kazimir Malevich\protect\footnotemark{} (right)%
  }%
  \label{fig:sliverart}
\end{figure}

In numerical simulations, slivers are unacceptable because they can degrade the accuracy of the finite element and finite volume approximations.
The standard state-of-the-art solution for eliminating slivers is to locally violate the Delaunay empty sphere condition.
If the Delaunay condition should be satisfied, which is particularly the case if using Voronoi finite volume methods fulfilling the maximum principle~\cite{Gaertner-0,Gaertner}, one encounters a notorious and frustrating feature of the current Delaunay triangulation tools: the creation of artificial slivers.
For example, triangulating vertices of a cubic lattice might result in a substantial magnitude of slivers (\cref{fig:tetgen:slivers}), the number and placement of which depend quite unpredictably on the rounding errors in the input data, creating although picturesque (\cref{fig:sliverart}) but unacceptable patterns even for the most simple input data.
Obtaining a Delaunay triangulation by means of a radical partition is more complicated but it allows for the elimination of artificial slivers without the violation of the Delaunay condition.

\footnotetext{Kazimir Malevich, ``Supremus~\textnumero~38'', public domain, via Wikimedia Commons.}

Note, however, that the same notation ``sliver'' is sometimes used for a badly shaped simplex which is not far from being flat but the computation of its circumcenter and circumradius is stable.
Such slivers are legitimate badly shaped Delaunay simplexes and should not be eliminated from the triangulation.
The difference between a correct but almost flat simplex and an incorrect triangulation of a Delaunay polyhedron face has to be clearly identified.
A natural identification is provided by the power diagram (radical partition) for the set of spheres, which goes back to René Descartes' work.
A \emph{radical partition} is defined as the set of convex polyhedra constructed by intersections of the half-spaces defined by the $d-1$-dimensional planes orthogonal to the segments connecting centers for all pairs of spheres (see \cref{sec:2} for details).
For two intersecting spheres, the radical plane always passes through their common set.
Hence, for a set of Delaunay spheres, the radical partition coincides with the Delaunay partition.

The equivalence of the radical and the Delaunay partitions provides a natural way to introduce a polyhedral approximation of a Delaunay partition for a disturbed point set.
Let us add a small perturbation to the vertex positions and consider an arbitrary simplicial Delaunay partition for this vertex set, also allowing slivers.
Evidently, this would  change the number of the Delaunay spheres as well as their centers and radii.
The centers and radii of the Delaunay spheres are, however, stable with respect to perturbations since the volume of each Delaunay polyhedron was strictly positive initially.
Hence, for each of the Delaunay spheres of the initial partition, we obtain a cluster of spheres approximating it.
A new sphere is the circumsphere of a certain new Delaunay simplex.
We can find the best fit for the center of each cluster using averaging of circumcenters with volumes of the corresponding Delaunay simplexes as weights.
In this case, evidently, a sliver would  provide a zero contribution to the averaged sphere parameters.
The most simple way to compute the average radius for a cluster of spheres is to use least squares for distances between the new center and the disturbed vertices.
A sliver might produce an isolated circumsphere which far from the cluster of spheres but such a sliver should be ignored when constructing the radical partition.
The elimination criterion is simple: if the circumsphere is unstable, then the almost flat simplex is too badly conditioned to contribute to the set of spheres.

\begin{figure}[t]%
  \centering{}%
  \includegraphics[width=0.48\textwidth]{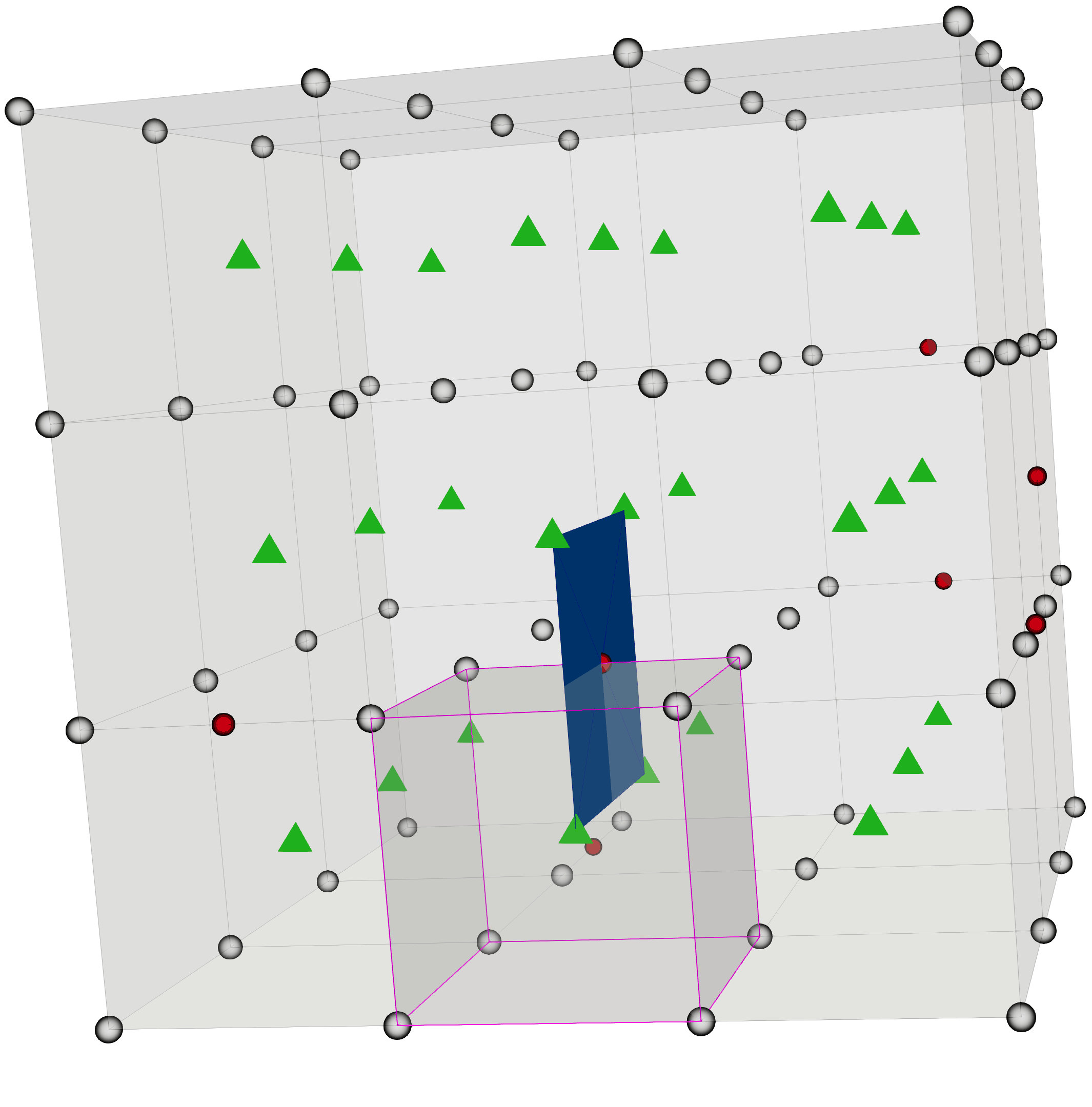}%
  \hfill{}%
  \includegraphics[width=0.48\textwidth]{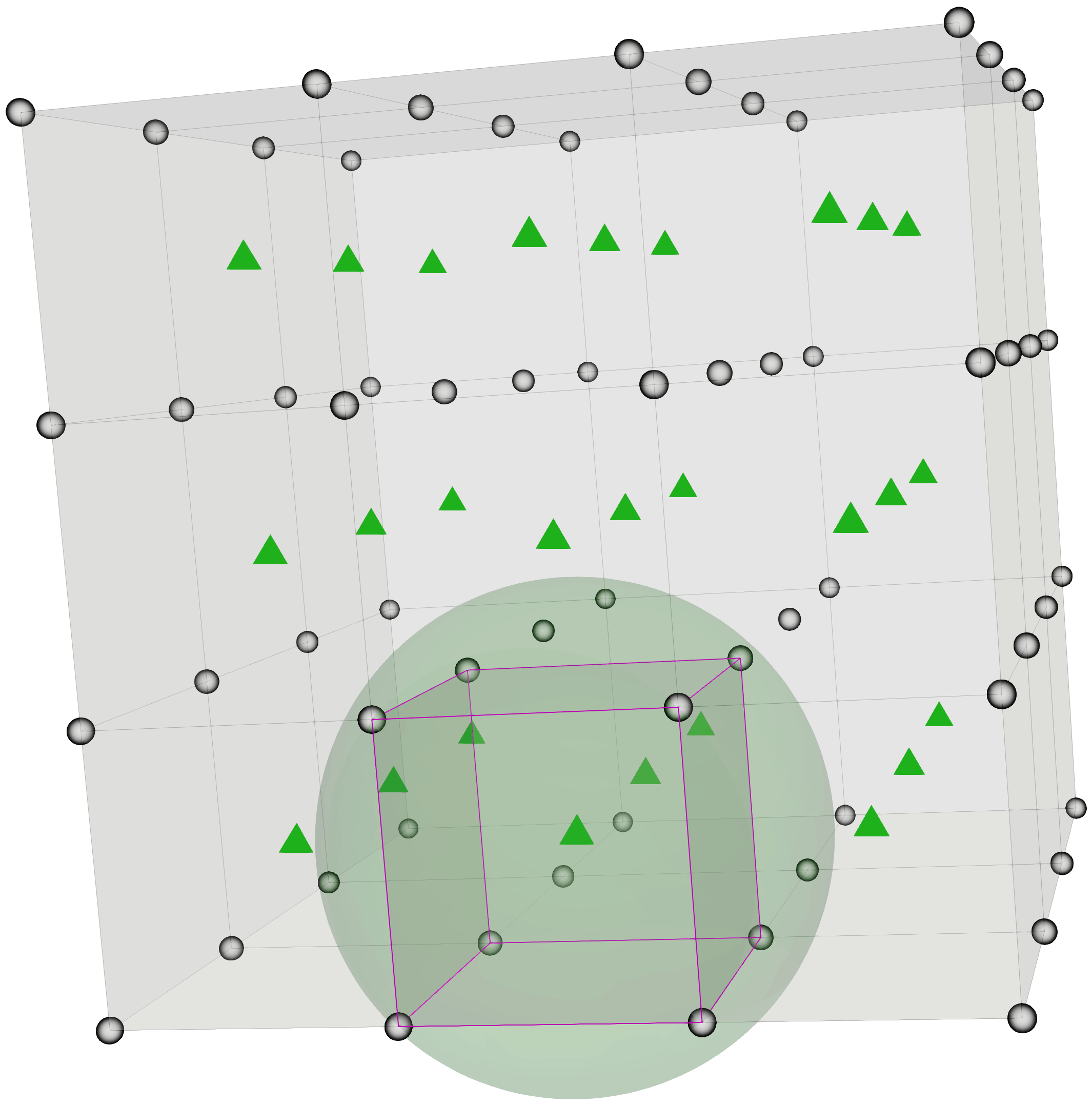}%
  \caption{Left: an artificial sliver (due to the floating point arithmetics
    rounding) creates a superfluous vertex of the radical partition.
    Right: the correct radical cell coinciding with the Delaunay cell}%
  \label{fig:tetgen:cubic}%
\end{figure}

\emph{TetGen}~\cite{TetGen} can be used efficiently to create a radical partition providing an option for constructing weighted Delaunay tetrahedralization.
In our case, vertices and weights are defined by the centers and the radii of the stable spheres (see details in \cref{sec:2}). 
A radical face is built by connecting orthocenters of the weighted Delaunay tetrahedra adjacent to an weighted Delaunay edge.
At this stage, \emph{TetGen} tends to produce artificial slivers, thus, adding superfluous vertices to radical faces (\cref{fig:tetgen:cubic}, left).
Such vertices and the tetrahedra that generate them can be simply ignored, resulting in a correct radical partition coinciding with the Delaunay partition of the original point set (\cref{fig:tetgen:cubic}, right).

As soon as a radical partition is computed, a stable triangulation of the resulting polyhedral partition can be found, which can be used as an approximate Delaunay triangulation.
Note, that the radical partition for the disturbed point set may result in another point set as a set of partition vertices.
It may produce a cluster of close vertices in place of an input vertex, potentially creating needle-type simplexes.
These clusters should also be glued together.

Delaunay triangulations have numerous applications in different fields~\cite{Aurenhammer}.
One of the very important properties for numerical simulations is the maximum principle for the discrete Laplacian on Delaunay meshes~\cite{Gaertner-0,Gaertner}.
They key ingredient for constructing truly Voronoi computational meshes without cutting them near the boundaries is the ability to control the location of Delaunay spheres in the key regions of the computational domain, in particular on and near the boundaries~\cite{Garanzha-2019-0,Garanzha-2019,vorocrust-2020}.
We propose a computational framework which, potentially, can serve as a controlling tool for the placement of Delaunay spheres and provide a top-to-bottom Delaunay-Voronoi construction when the set of spheres generates the set of seeds for Delaunay meshing.

\section{Power Diagram and the Lifting Procedure}%
\label{sec:2}

The idea of lifting goes back to the works of G.F.~Voronoi~\cite{Voronoi-1908}, who showed that a Delaunay triangulation in $\R^d$ is the projection of faces of a convex polyhedron $P \in R^{d+1}$ inscribed into a circular paraboloid $\Pi$.
The convex body $P^*$, constructed as the intersection of the upper half-spaces for the tangent planes to $\Pi$ at the vertices of $P$, is called the \emph{Voronoi generatrice}.
The projection of its faces onto $\R^d$ defines the \emph{Voronoi diagram}.
A more general lifting concept~\cite{Edelsbrunner-1,Edelsbrunner-2} is based on the construction of a pair $P$, $P^*$ of convex polyhedra which satisfy the polarity relation~\cite{Alexandrov-convex} with respect to the paraboloid $\Pi$.

Consider the system of balls $\mathcal{B} = \qty{\, B_1, \dotsc, B_n \,}$ defined by  centers $\vb*{c}_i \in \R^d$ and radii $R_i \geq 0$,
the lifted point system $\mathcal{E}_l = \qty{\, \vb*{p}_1, \dotsc, \vb*{p}_n \,} $ in $\R^{d+1}$, where
\[
   \vb*{p_i}^T
   = \qty(\vb*{c}_i^T, \frac{1}{2}\qty(\abs{\vb*{c}_i}^2 - R_i^2) )
   = \qty(\vb*{c}_i^T,  h_i)
  ,
\]
and the lower convex envelope of $\mathcal{E}_l$ defined by the convex function $x_{d+1} = v(x_1, \dots, x_d)$.
The Legendre-Young-Fenchel dual~\cite{Fenchel} of $v$ is denoted by $v^*$.
To be mathematically precise, the function $v(\vb*{x})$ is equal to $+ \infty$ outside the convex envelope $\conv (\vb*{c}_1, \dotsc, \vb*{c}_n)$ and its epigraph (supergraph) is a closed set.
The dual function $v^*(\vb*{x})$ is defined everywhere and its graph contains unbounded faces.
In the following, we will consider constrained problems where unbounded faces are essentially excluded from the problem setting.

The projection of faces of the graph of $v$ defines a weighted Delaunay triangulation $\mathcal{W}$ in $\R^d$~\cite{Edelsbrunner-2}.
The vertices of the graph of $v$ are pairs $(\vb*{c}_i, h_i)$ and the vertices of the weighted Delaunay triangulation are $\vb*{c}_i$.
$T_k$ denotes the $k$-th weighted Delaunay polyhedron.

\begin{figure}[t]%
  \raisebox{0cm}{\includegraphics[width=0.46\textwidth]{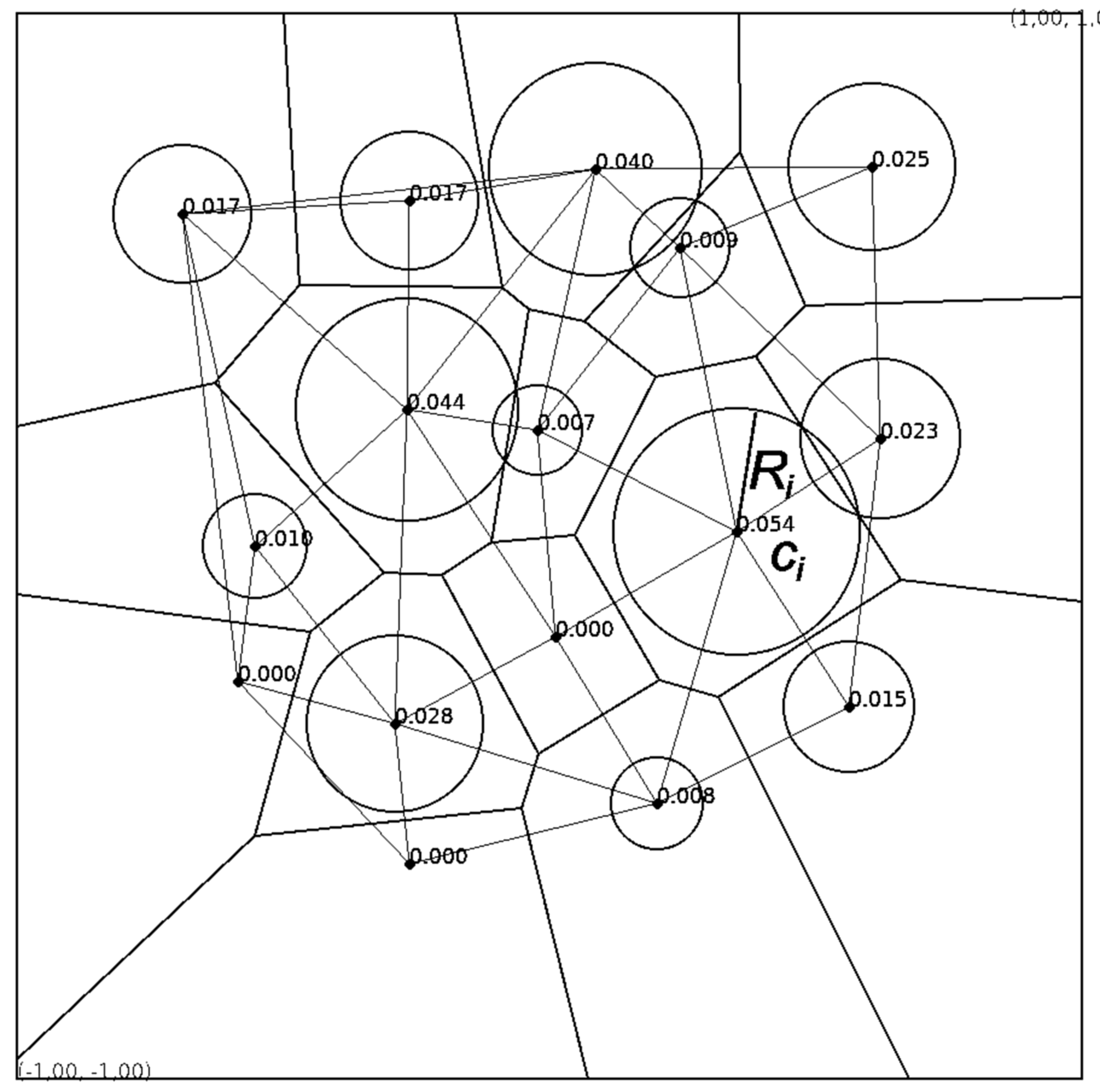}}%
  \hfill{}%
  \includegraphics[width=0.51\textwidth]{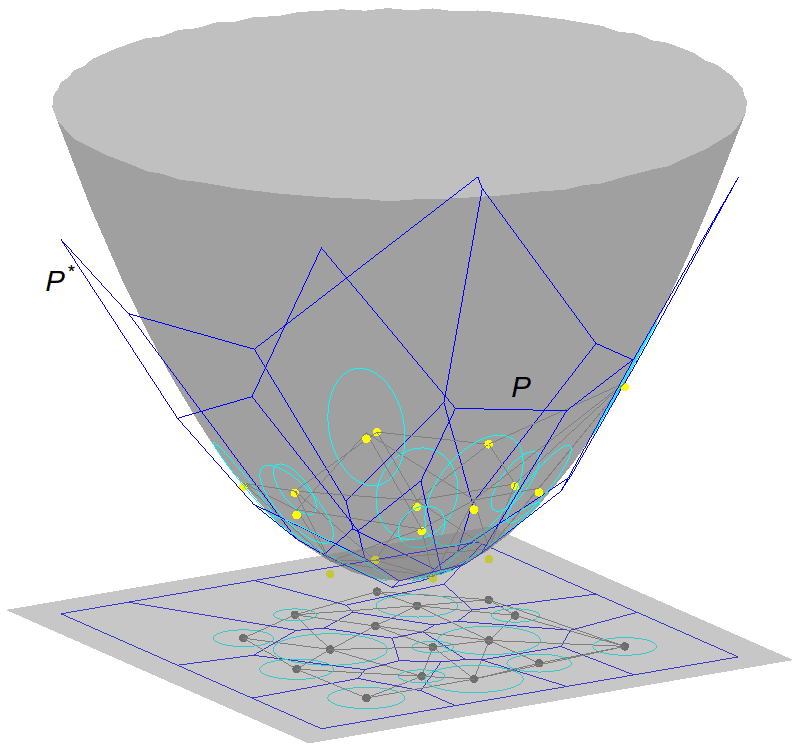}%
  \caption{A weighted Delaunay triangulation, the power diagram, and dual polyhedra from lifting.
    The figure is constructed with the help of \texttt{detri2} by Hang Si}%
  \label{fig:lifting}%
\end{figure}

The projection of the faces of the graph of $v^*$ defines a radical partition (power diagram) $\mathcal{R}$ for a system of balls in $\R^d$~\cite{Edelsbrunner-2} (\cref{fig:lifting}).
The projection of the $k$-th vertex of the graph of $v^*$ onto $x_{d+1} = 0$ is denoted by $\vb*{v}_k$.
This point is dual to $T_k \in \mathcal{W}$, while the vertex $\vb*{c}_i$ is dual to the cell $D_i$ of $\mathcal{R}$.

\section{Ball Movement as a Transformation of Dual Polyhedra}

Our objective is to move and scale the balls $\mathcal{B}$ in such a way that all vertices of the graph of $v^*$ converge to the surface of the paraboloid $x_{d+1} = \Pi(\vb*{x}) = \qty(x_1^2 + \cdots + x_d^2)/2$.
It means that the projection of the graph of $v^*$ will eventually converge to the Delaunay partition.
Note, that the number of vertices $\vb*{v}_k$ may vary during the ball movement.
Further, it is convenient to use the notation $(\vb*{x}^T, x_{d+1})$ for an arbitrary point in $\R^{d+1}$, where $\vb*{x}^T = \qty{ x_1, \dots, x_d }$.

For the set of balls $\mathcal{B}$, we build the primal and the dual polyhedra $P$ and $P^*$ defined by the convex piecewise-linear functions $v(\vb*{x})$ and $v^*(\vb*{x})$, respectively.
According to the polarity relation with respect to the circular paraboloid~\cite{Alexandrov-convex}, the vertex $(\vb*{c}_i, h_i)$ defines the dual plane of the face of the graph of $v^*(\vb*{x})$,
\[
   \vb*{c}_i^T \vb*{x} = h_i + x_{d+1}
  .
\]
The vertex $(\vb*{v}_k, z_k)$ of the graph of $v^*(\vb*{x})$ is the intersection of at least $d+1$ such planes:
\begin{equation}
\vb*{v}_k^T (\vb*{c}_i - \vb*{c}_p) = h_i - h_p
  ,
  \label{2028:eq.dual-vertex}
\end{equation}
where $i \neq p$ are indices of all planes intersecting in $(\vb*{v}_k, z_k)$.
Hence,
\[
  z_k
  = \vb*{v}_k^T \vb*{c}_i - h_i
  = \vb*{v}_k^T \vb*{c}_i - \frac{1}{2} \vb*{c}_i^2 + \frac{1}{2} R_i^2
  = - \frac{1}{2} \qty( \qty|\vb*{c}_i - \vb*{v}_k|^2 - R_i^2) + \frac{1}{2} \vb*{v}_k^2
  .
\]
In another words,
\begin{equation}
   z_k - \Pi(\vb*{v}_k) = -\frac{1}{2} \tau_i(\vb*{v}_k)
  \label{2028:ed.distance-paraboloid}
  ,
\end{equation}
where
\[
   \tau_i(\vb*{y}) = \qty|\vb*{c}_i - \vb*{y}|^2 - R_i^2
\]
is the power of the point $\vb*{y}$ with respect to the ball $B_i$.
Hence, the vertical distance of the vertex $(\vb*{v}_k, z_k)$ from the paraboloid $\Pi$ is fully defined by the value of the power.
Another interpretation of the \cref{2028:eq.dual-vertex} is that for a vertex $\vb*{v}_k$ dual to the weighted Delaunay polyhedron $T_k$ the equality
\[
  \tau_i(\vb*{v}_k) = \tau_p(\vb*{v}_k)
\]
is satisfied for all vertices of $T_k$.
It means that, in our setting, it is possible to omit the indices and use the notation $\tau(\vb*{v}_k)$ for the value of the power.

\begin{figure}%
  \centering{}%
  \subfloat[\label{fig:radical:split:a}]{%
    \includegraphics[width=0.40\textwidth]{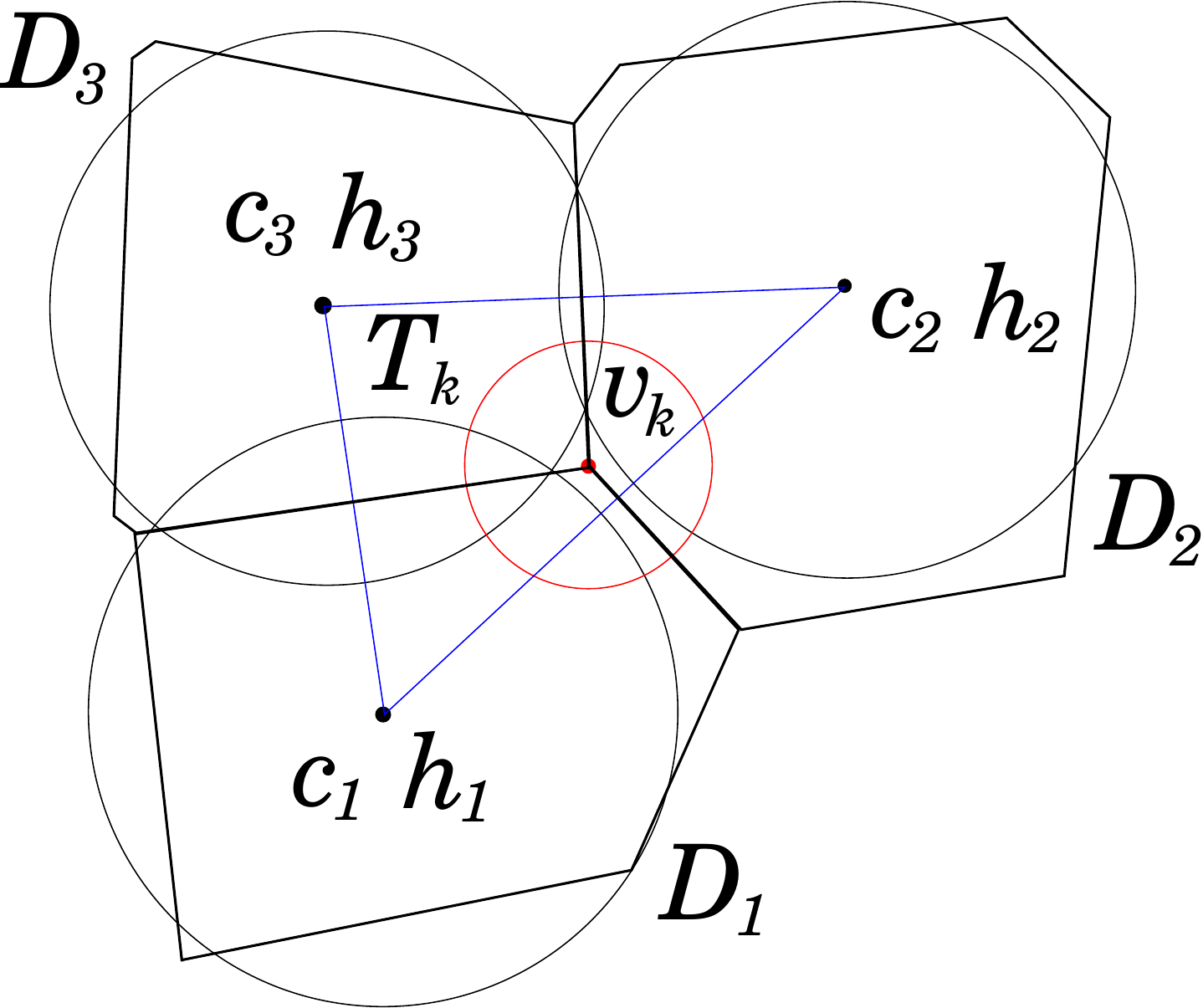}%
  }%
  \hspace{0.1\textwidth}%
  \subfloat[\label{fig:radical:split:b}]{%
    \includegraphics[width=0.36\textwidth]{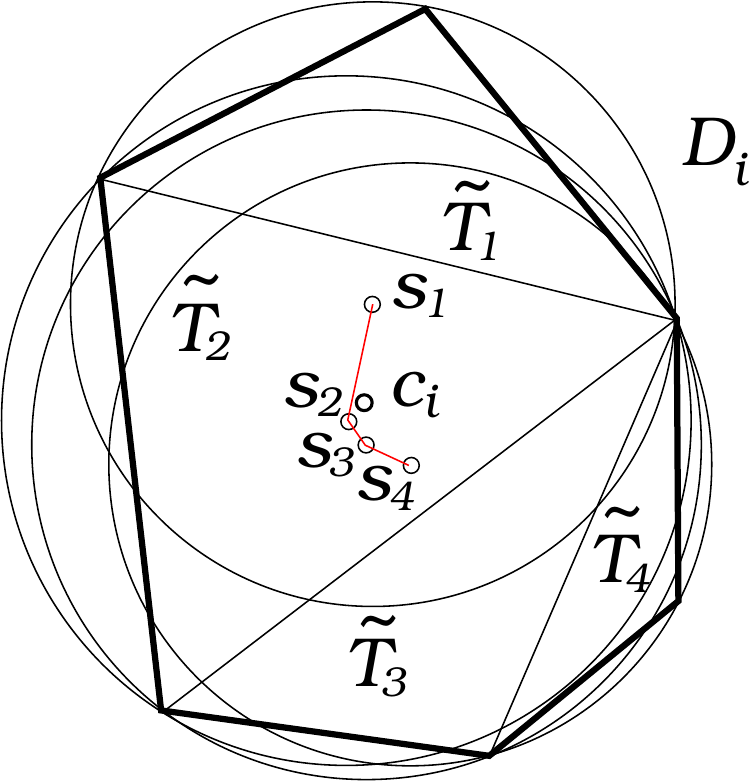}%
  }%
  \caption{%
    \subref{fig:radical:split:a}
    A weighted Delaunay triangle $T_k$, the dual vertex $\vb*{v}_k$ and orthocircle,
    \subref{fig:radical:split:b}
    a radical cell is split into  $4$ Delaunay triangles.
    Center $\vb*{c}_i$ is approximated by a cloud of $4$ Delaunay circumcenters}%
  \label{fig:radical:split}%
\end{figure}

Relation \cref{2028:eq.dual-vertex} implies that the gradient of $v^*(\vb*{x})$ at the $i$-th face of its graph is equal to $\vb*{c}_i$.
The dual statement is true as well: the gradient of $v(\vb*{x})$ at the $k$-th face of its graph is equal to $\vb*{v}_k$.
For any convex polyhedron $T_k$ there are $d$ linear independent vectors $\vb*{c}_i - \vb*{c}_p$.
In two dimensions, in the simplest case, $T_k$ is a triangle and this set is $\qty{\vb*{c_2} - \vb*{c}_1, \vb*{c}_3 - \vb*{c}_1}$ (\cref{fig:radical:split:a}).
When $\tau(\vb*{v}_k) > 0$, one can associate with $\vb*{v}_k$ a sphere with the radius $\sqrt{\tau(\vb*{v}_k)}$ which is called the \emph{orthosphere}.
In our numerical experiments in \cref{sec:numerical}, we draw artificial orthospheres with radii defined by $\sqrt{\qty|\tau(\vb*{v}_k)|}$.
As shown below, these spheres visualize the deviation of the radical partition from the Delaunay partition.

Let $\tilde{v}^*$ be the projected version of the function $v^*$ constructed as follows: consider the set of vertices $\qty{\, (\vb*{v}_j, z_j) \,}$ of the graph of $v^*$ and project them on the paraboloid $\Pi$ by setting
\[
  \tilde{z}_j = \frac{1}{2} \qty|\vb*{v}_j|^2
  \qq{and}
  \tilde{\vb*{v}}_j = \vb*{v}_j
  .
\]
As shown above, this projection changes the vertical component $z_j$ by $\tau(\vb*{v}_j)/2$.

Computing the lower convex envelope of the point system $\qty{\, (\tilde{\vb*{v}}_k, \tilde{z}_k) \,}$, $k = 1, \dotsc, n_v$, results in the graph of the function $\tilde{v}*$.

\section{Dirichlet Functional for Power Function}

In order to measure the closeness of the current radical partition to the Delaunay partition, we consider the Dirichlet functional for the difference of $v^*$ and $\tilde{v}^*$,
\begin{equation}\label{2028:eq.Dirichlet}
   F(X) = \frac{1}{2} \int_\Omega \qty|\nabla v^* - \nabla \tilde{v}^*|^2 \dd{\vb*{x}}
  .
\end{equation}
Here, $X$ is the vector of unknowns, consisting of $\vb*{c}_i$ and $R_i$, and $\Omega$ is the bounded definition domain of the function $\tilde{v}^*$.
From \cref{2028:ed.distance-paraboloid},
\[
   F(X) = \frac{1}{8} \int_\Omega \qty|\nabla \tau(\vb*{x})|^2 \dd{\vb*{x}}
  .
\]
Let $\tau(\vb*{x})$ denote a piecewise linear function which coincides with $\tau(\vb*{v}_k)$ at $\vb*{v}_k$ and is linear on each cell $\tilde{T}_j$ of the auxiliary Delaunay partition.
Assuming that each face of the graph of $v^*$ is projected onto the paraboloid independently of other faces, $F(X)$ can be rewritten as
\begin{equation}
  \label{2028:eq.Dirichlet-final}
  F_I(X) = \frac{1}{2} \sum_i \sum_{\tilde{T}_j \in D_i}
  \qty|\vb*{c}_i - \vb*{s}_j|^2 \vol {\tilde{T}_j} ,
\end{equation}
where $\vb*{s}_j$ is the circumcenter of the Delaunay simplex $\tilde{T}_j$.
The above equality is the obvious consequence of the fact that
\[
  \nabla v^* \, \vert_{D_i} = \vb*{c}_i
  \qc
  \nabla \tilde{v}^* \, \vert_{\tilde{T}_j \in D_i} = \vb*{s}_j
  .
\]
To clarify this formula, consider a Delaunay simplex $\tilde{T}_j$ with vertices $\vb*{v}_1, \dotsc, \vb*{v}_{d+1}$.
The gradient $\vb*{g}_j$ of  $\tilde{v}^*$ is defined by
\[
   {\qty(\vb*{v}_m - \vb*{v}_l)}^T \vb*{g}_j 
  = \frac{1}{2} \qty(\vb*{v}_m^2 - \vb*{v}_l^2)
  \qc m, l \leq d+1
  ,
\]
or
\[
  \qty|\vb*{g}_j - \vb*{v}_m|^2 = \qty|\vb*{g}_j - \vb*{v}_l|^2
  ,
\]
which is precisely the set of equations for the circumcenter $\vb*{s}_j$ of $\tilde{T}_j$.
Note, that in general the functional $F_I(X)$ does not coincide with $F(X)$, since the independent projection of faces onto the paraboloid might result in an inscribed polyhedron which might loose the convexity on the boundary of projected faces.
However, upon convergence, $F_I$ and $F$ coincide.

The equality $F_I(X) = 0$ implies that for each radical polyhedron $D_i$ its dual vertex $\vb*{c}_i$ coincides with all circumcenters of the Delaunay triangulation of its set of vertices, meaning that $D_i$ is a Delaunay polyhedron.

\smallskip{}

Consider the following algorithm:

\begin{itemize}[itemsep=0ex]
  \item[]For $n = 0, 1, \dots$
  \item Given the set of balls $\mathcal{B}^n$, compute the primal and dual functions $v^n$ and ${v^n}^*$ using lifting on the paraboloid.
  These functions define a weighted Delaunay triangulation $\mathcal{W}^n$ and a radical partition $\mathcal{R}^n$, respectively.

  \item Compute the projected function ${\tilde{v}^n} {}^*$.
  This function defines a Delaunay triangulation $\tilde{\mathcal{T}}^n$ of the set of vertices of radical partition $\mathcal{R}^n$.

  \item Do an approximate gradient descent minimization step for the Dirichlet functional $F_I(X)$ to obtain the set $\mathcal{B}^{n+1}$.

  \item[] Repeat until convergence.
\end{itemize}

In this algorithm, the sequence of weighted Delaunay triangulations $\mathcal{W}^n$ converges to a Voronoi diagram $\mathcal{V}$ while the sequence of radical partitions $\mathcal{R}^n$ converges to a Delaunay partition $\mathcal{T}$.
Note that the limit Delaunay triangulation $\tilde{\cal T}$ consists of Delaunay simplexes, while $\mathcal{T}$ consists of polyhedral Delaunay cells for the same point set.
It is not necessary to build a consistent Delaunay triangulation at each step of the algorithm.
If needed, it can be built only once using the final Delaunay partition.

For a simplified version of the algorithm avoiding the cumbersome computation of the exact gradient, the centers $\vb*{c}_i$ are computed using the local minimization of the functional \cref{2028:eq.Dirichlet-final} considering it as a quadratic function of $\vb*{c}_i$,
\[
   \vb*{c}_i^{new} = \qty(\sum_{\tilde{T}_j \in D_i}  \vb*{s}_j \vol {\tilde{T}_j} )
  ~\Big/ \sum_{\tilde{T}_j \in D_i}  \vol {\tilde{T}_j}
  ,
\]
while $R_i$ is computed using the simple least squares approximation
\begin{equation}\label{2028:newR}
   R_i^{new} = \frac{1}{\sqrt{M}} {\qty(\sum\limits_{m=1}^M \abs{\vb*{c}_i^{new} - \vb*{v}_m}^2)}^\frac{1}{2}
  ,
\end{equation}
where $\vb*{v}_1, \dotsc, \vb*{v}_M$ are vertices of $D_i$.
Note, that \cref{2028:newR} is equivalent to
\[
   \sum\limits_{m=1}^M \tau(\vb*{v}_m) = 0
  ,
\]
which, in turn, is the necessary minimum condition with respect to $R_i$ for the local functional
\[
  \sum\limits_{m=1}^M \tau^2(\vb*{v}_m)
  ,
\]
provided that the centers are fixed.
New positions and radii are computed using a relaxation parameter $0 < \theta < 1$,
\begin{align*}
   \vb*{c}_i^{n+1} &= \vb*{c}_i^n (1 - \theta) + \vb*{c}_i^{new} \theta,
  \\
  \ R_i^{n+1} &= R_i^n (1 - \theta) + R_i^{new} \theta
  .
\end{align*}
This heuristic algorithm is quite efficient for initial iterations.
Eventually, it slows down and should be replaced by the gradient descent technique in order to converge to the exact solution.

A constrained problem where some of the balls are fully or partially fixed can be considered as well: ball centers are allowed to move along a prescribed manifold, radii constraints are added, etc.

Currently, there is no existence theorem for this problem.
Moreover, an overdetermined constrained problem can be easily defined such that the set of Delaunay spheres cannot be constructed by means of admissible operations without allowing the introduction of new spheres.

\section{Numerical Experiments}%
\label{sec:numerical}

For a relatively simple test setting, we consider a square with a circle inside and cover the all internal and boundary curves by the protecting circles.
The centers of the circles for the outer boundary are fixed.
The points if their intersection, which are outside of the domain, essentially define a set of additional Delaunay vertices and, thus, the outer boundary by the Voronoi edges.
For the inner boundary (the inner circle), we fix both the position of the circle centers and the radii.
Additional two lines of circles are added, defining the initial layer of quadrilateral Delaunay cells near the internal boundary.
The centers of these additional circles are fixed but the radii are allowed to change.
We allow a small change of the radii on the outer boundary as well.
Since the outer boundary is essentially already covered by Delaunay cells, we actually set the zero boundary conditions for the power function $\tau(\vb*{x})$.
Inside the domain, a set of circles placed regularly is given and a random perturbation is added to the centers and the radii.
All circles, the centers of which are outside of the domain or inside the protecting circles are removed.

\begin{figure}[H]%
  \includegraphics[width=0.48\textwidth]{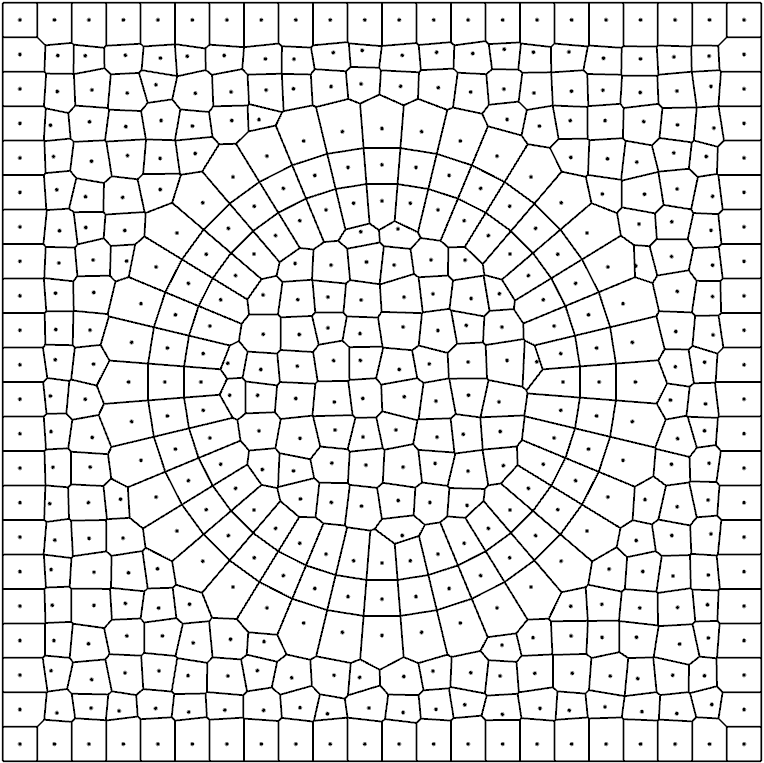}%
  \hfill{}%
  \includegraphics[width=0.48\textwidth]{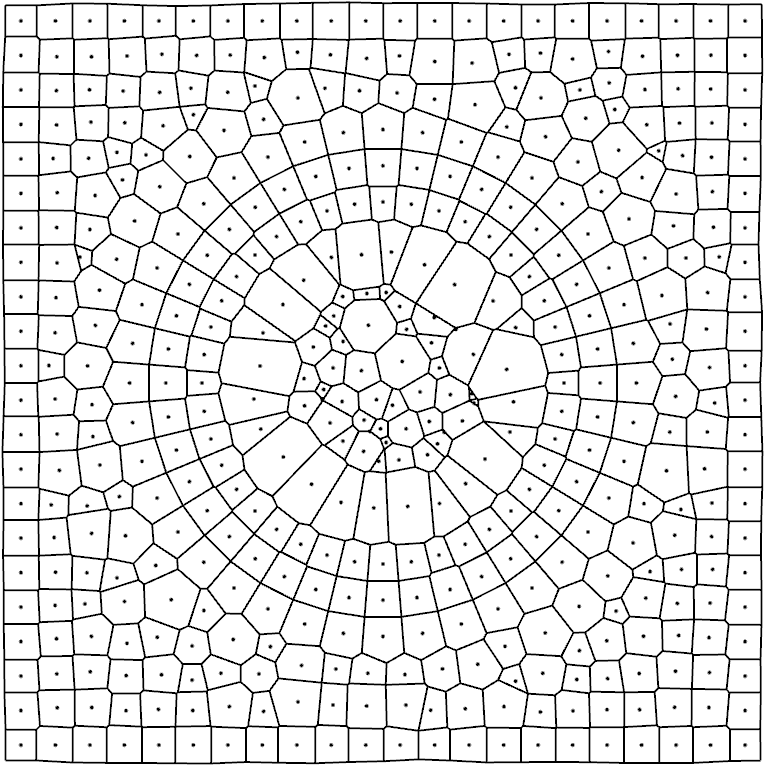}%
  \caption{Initial radical partition $\mathcal{R}_0$ and the final Delaunay partition $\mathcal{T}$}%
  \label{fig:Delone}%
\end{figure}

\Cref{fig:Delone} shows the initial radical partition and the final Delaunay partition.
It can be observed that the constraints are satisfied due to the appearance of larger cells near the fixed circles which, in turn, compresses the free circles almost to zero.

\begin{figure}[H]%
  \includegraphics[width=0.48\textwidth]{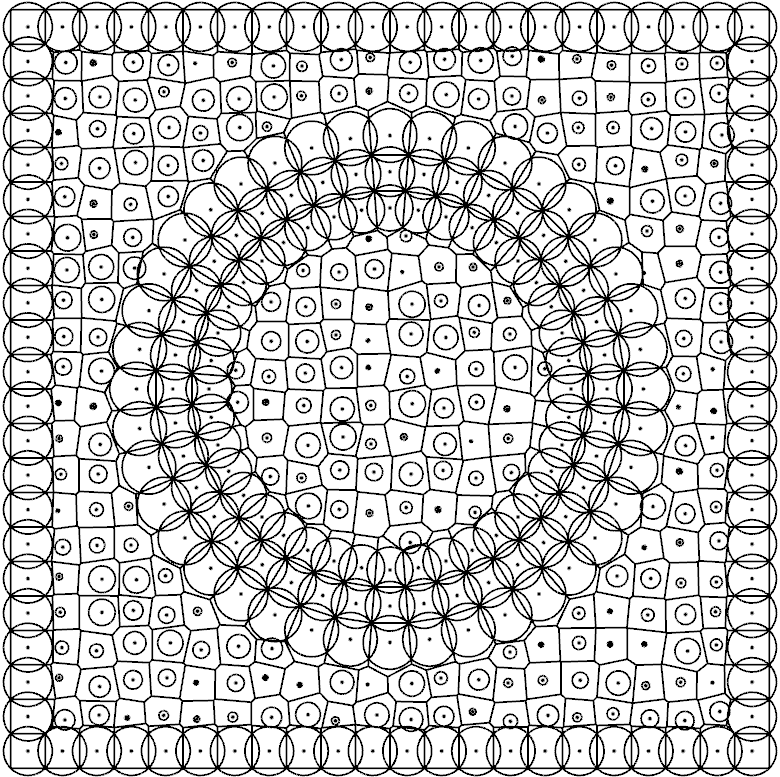}%
  \hfill{}%
  \includegraphics[width=0.48\textwidth]{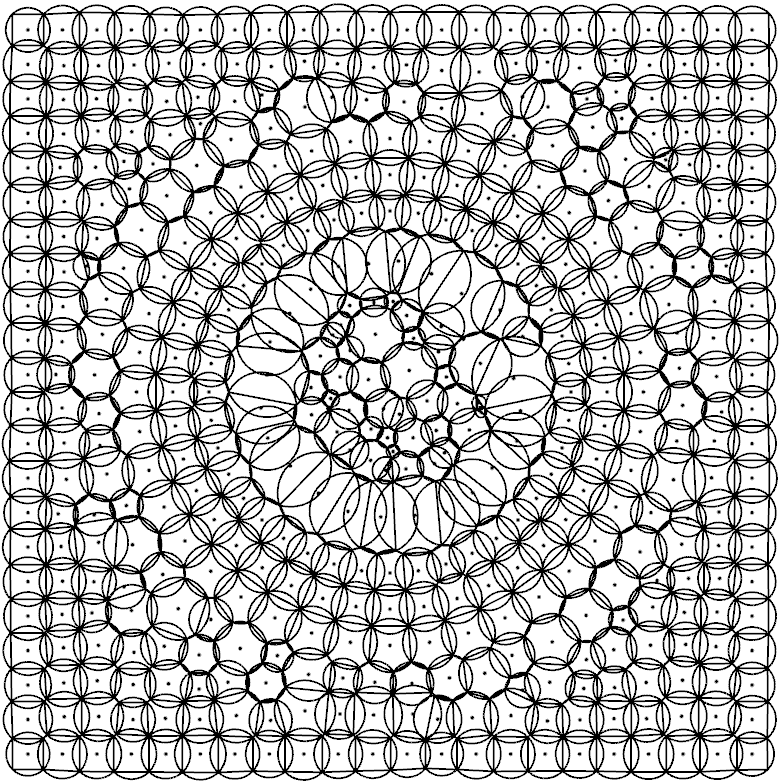} %
  \caption{The initial radical partition $\mathcal{R}_0$ with the initial set of circles and the final Delaunay partition $\mathcal{T}$ with the Delaunay circles}%
  \label{fig:Delone-circ}%
\end{figure}

\newpage{}%
Note, that the regular pattern of Delaunay and Voronoi cells near the internal boundary is disturbed by the appearance of small Delaunay edges.
In the current version of the algorithm, there is no mechanism to eliminate these.
\Cref{fig:Delone-circ} shows radical partitions and the circles $B_i$ that generate them.

\begin{figure}[H]
  \includegraphics[width=0.48\textwidth]{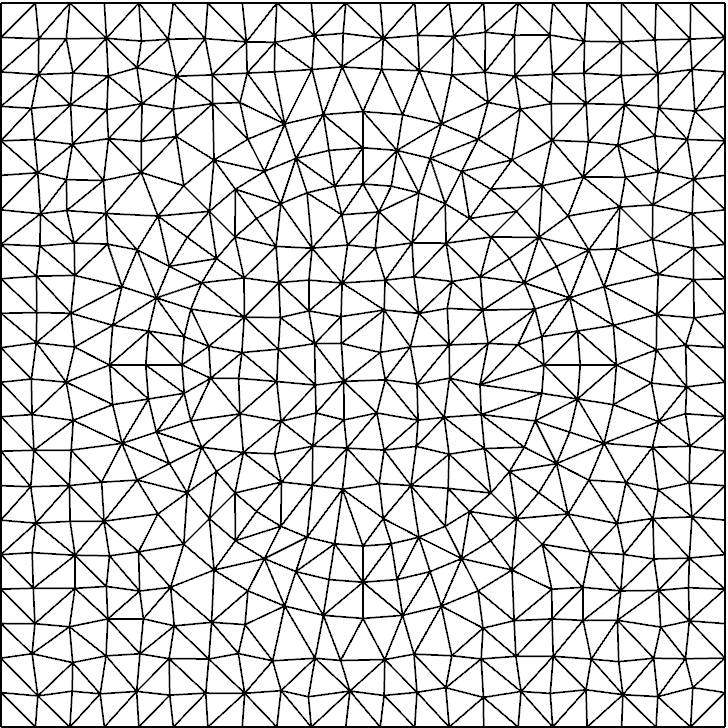}%
  \hfill{}%
  \includegraphics[width=0.48\textwidth]{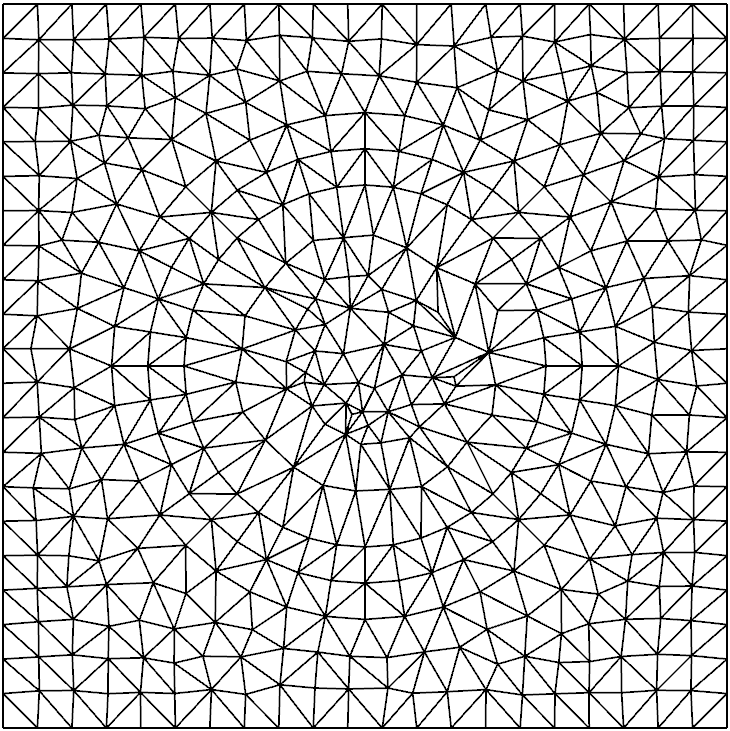} %
  \caption{The initial weighted Delaunay mesh $\mathcal{W}_0$ and the final Voronoi triangulation}%
  \label{fig:vorotri}%
\end{figure}

Initial weighted Delaunay triangulation $\mathcal{W}_0$ and the final Voronoi mesh are shown in \cref{fig:vorotri}.
Again, we observe quality problems in final Voronoi triangulation.
Note, that the two layers of Voronoi cells near the internal circle are, in fact, quadrilateral layers which are split into triangles for the presentation.

\begin{figure}[H]%
  \centering{}%
  \includegraphics[width=0.46\textwidth]{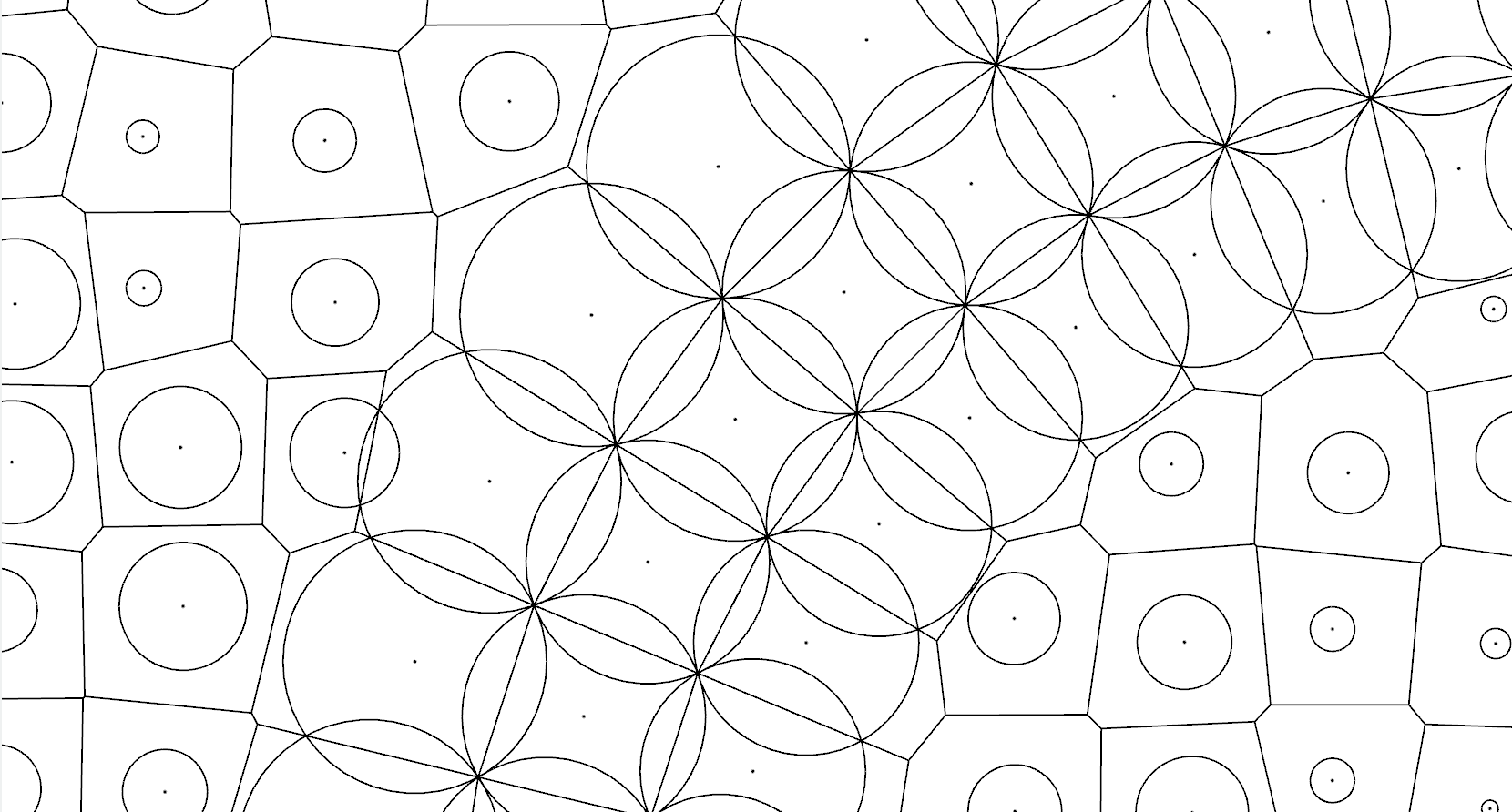}%
  \hspace{0.05\textwidth}%
  \includegraphics[width=0.46\textwidth]{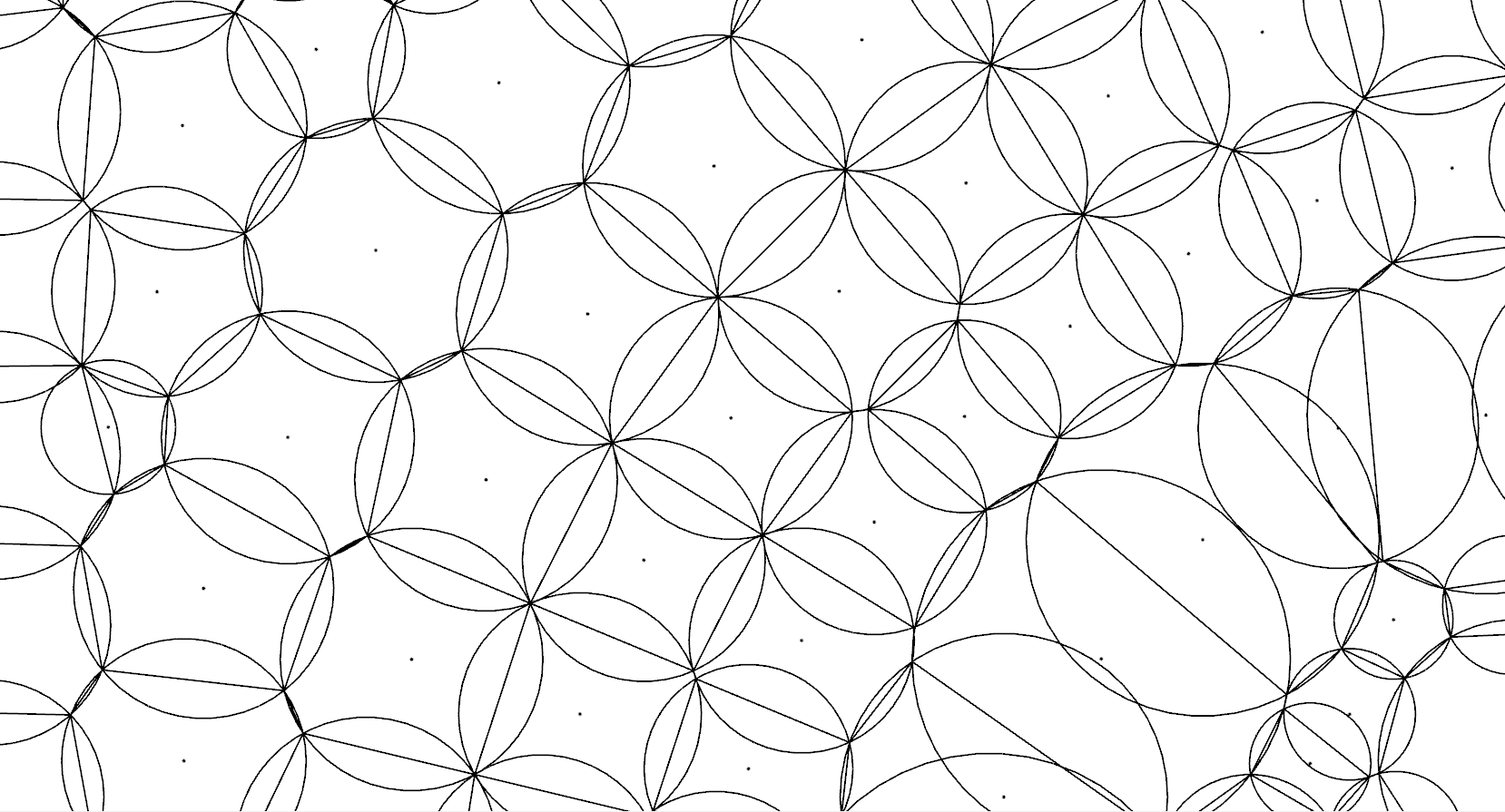}%
  \\[2ex]
  \includegraphics[width=0.46\textwidth]{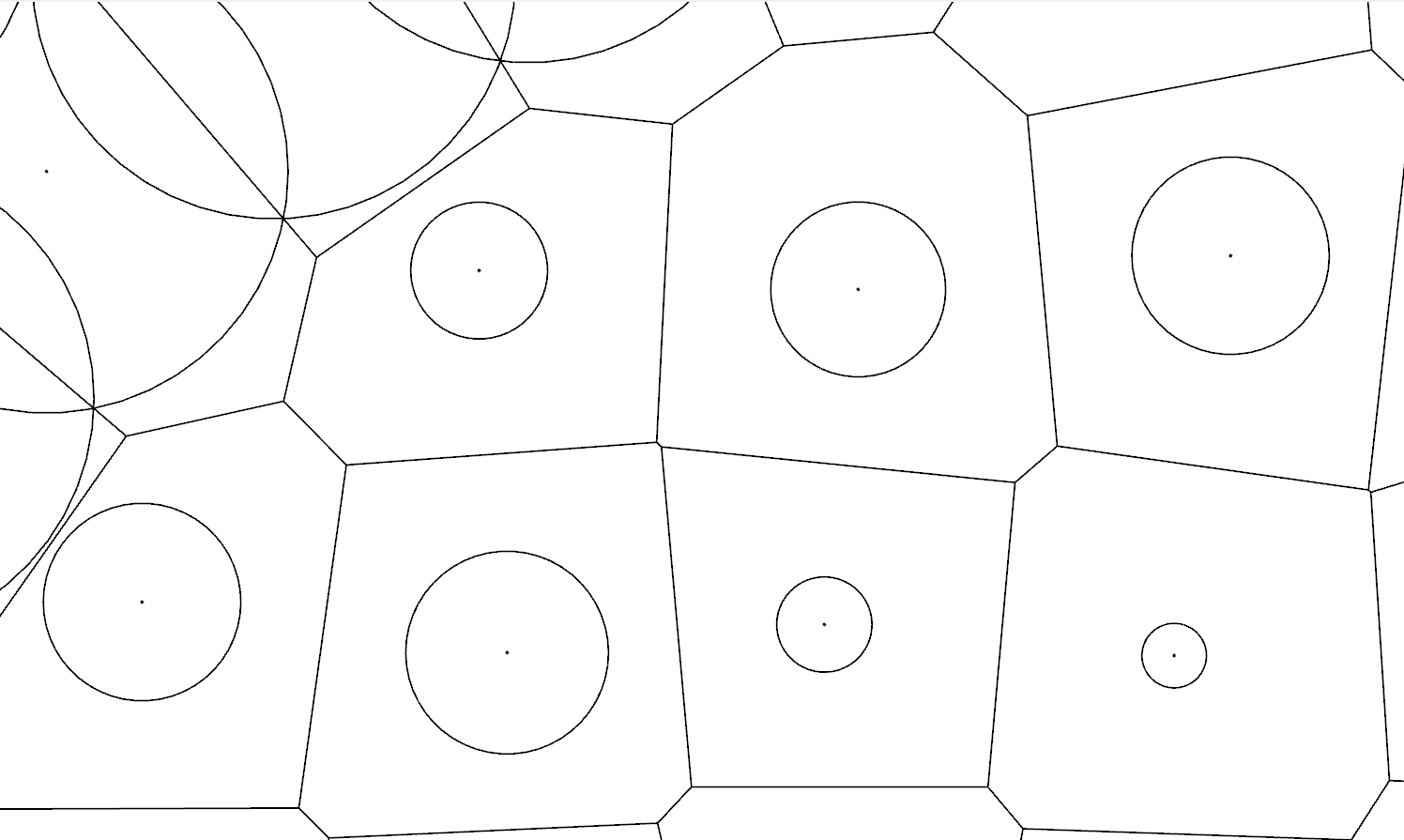}%
  \hspace{0.05\textwidth}%
  \includegraphics[width=0.46\textwidth]{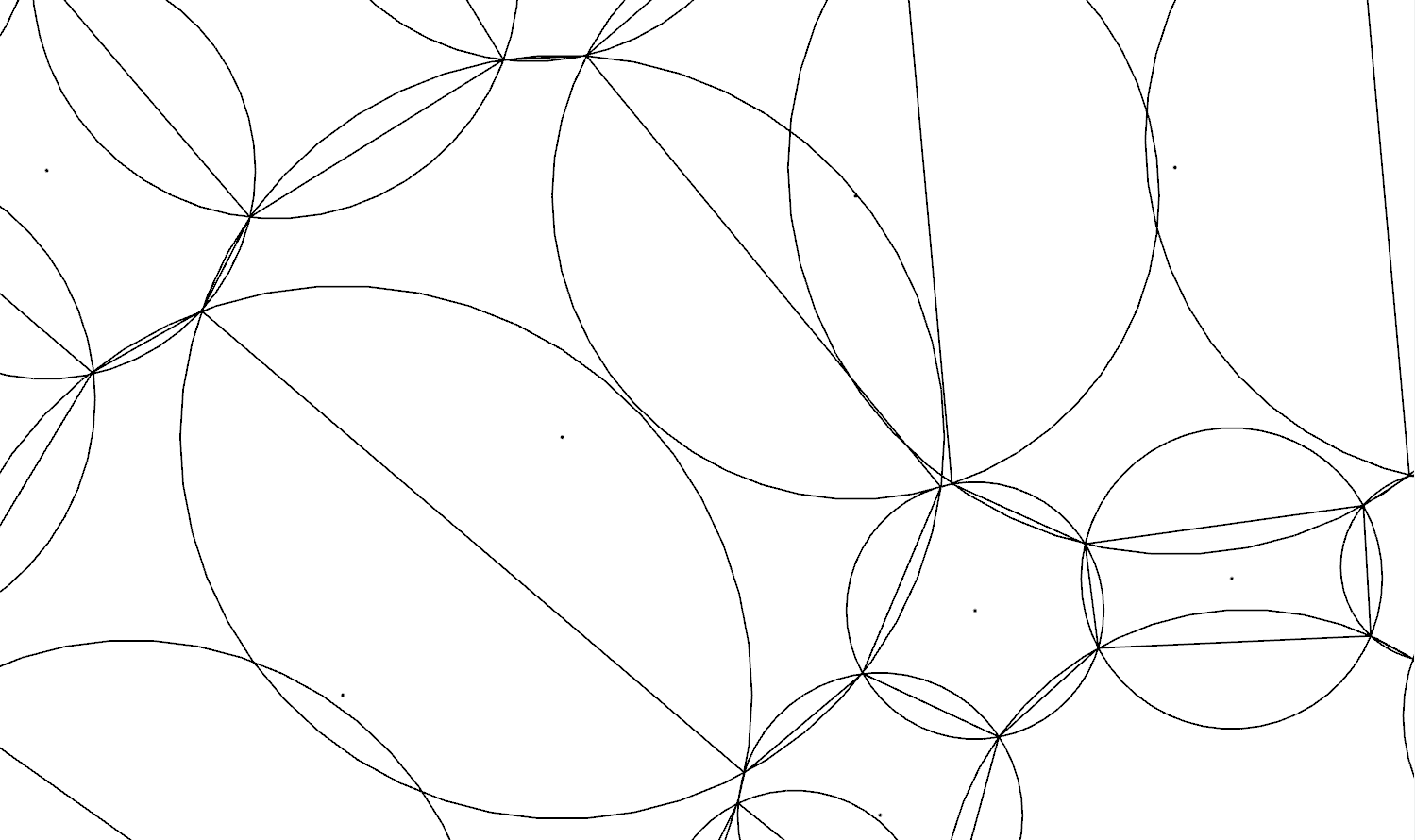}
  \caption{Fragments of the initial radical partition (left)
  and the corresponding fragments of the final Delaunay partition (right)}%
  \label{fig:Delone:ff}%
\end{figure}

Enlarged fragments of the meshes in \cref{fig:Delone:ff} allow to see clearly how the algorithm works.

\begin{figure}[H]%
  \centering{}%
   \subfloat[\label{fig:Delone-power:a}]{%
    \includegraphics[width=0.46\textwidth]{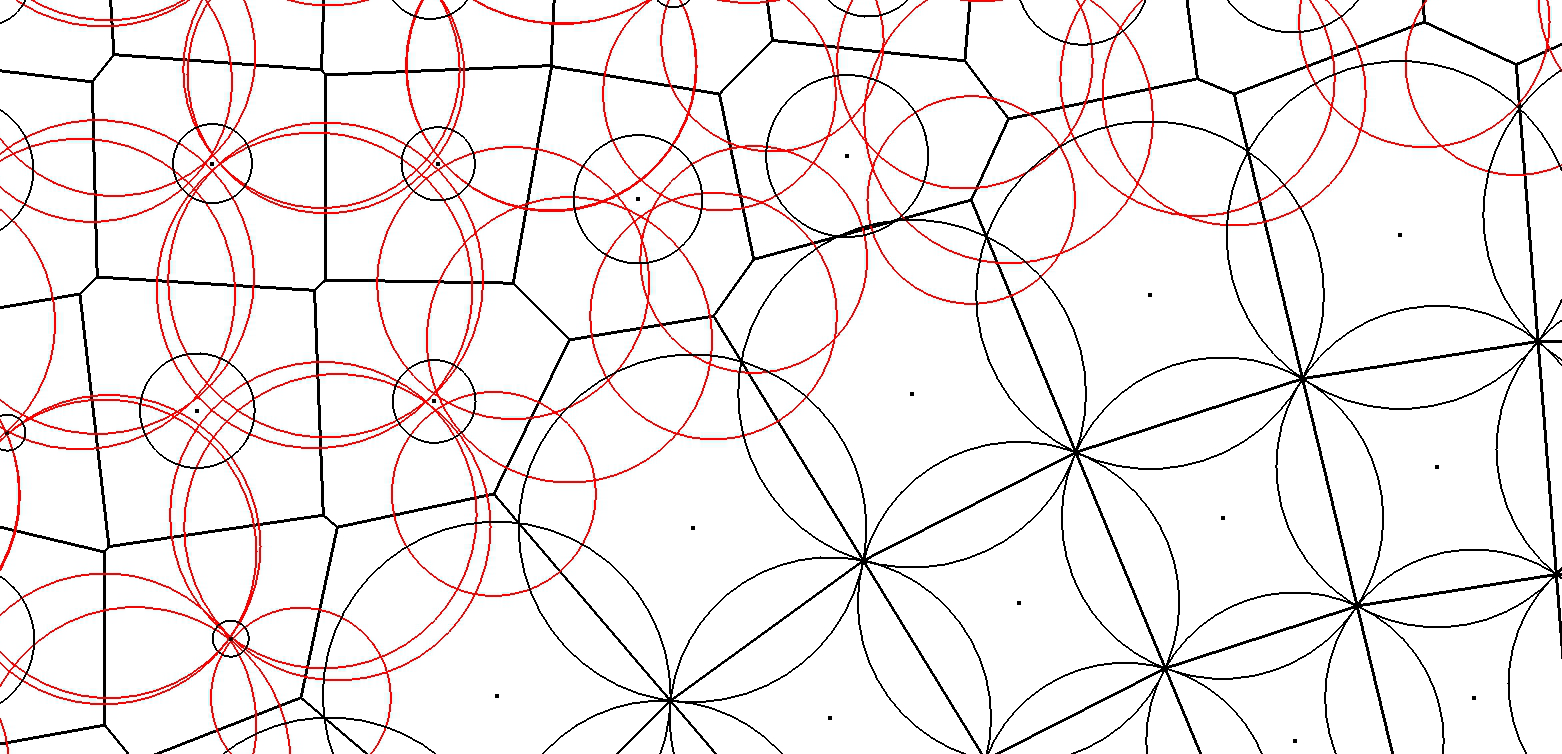}%
  }%
  \hspace{0.05\textwidth}%
  \subfloat[\label{fig:Delone-power:b}]{%
    \includegraphics[width=0.46\textwidth]{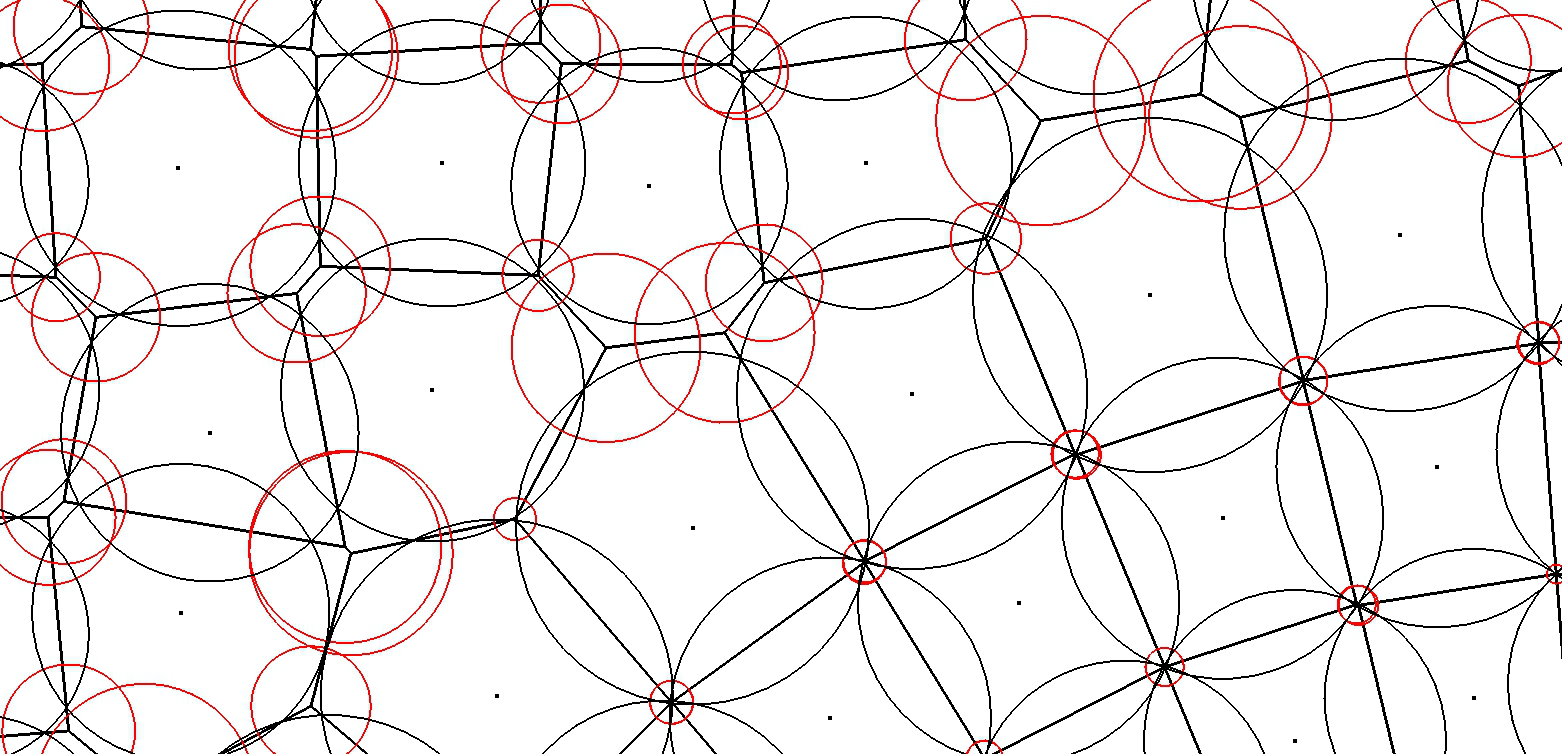}%
  }%
  \\[1ex]%
  \subfloat[\label{fig:Delone-power:c}]{%
    \includegraphics[width=0.46\textwidth]{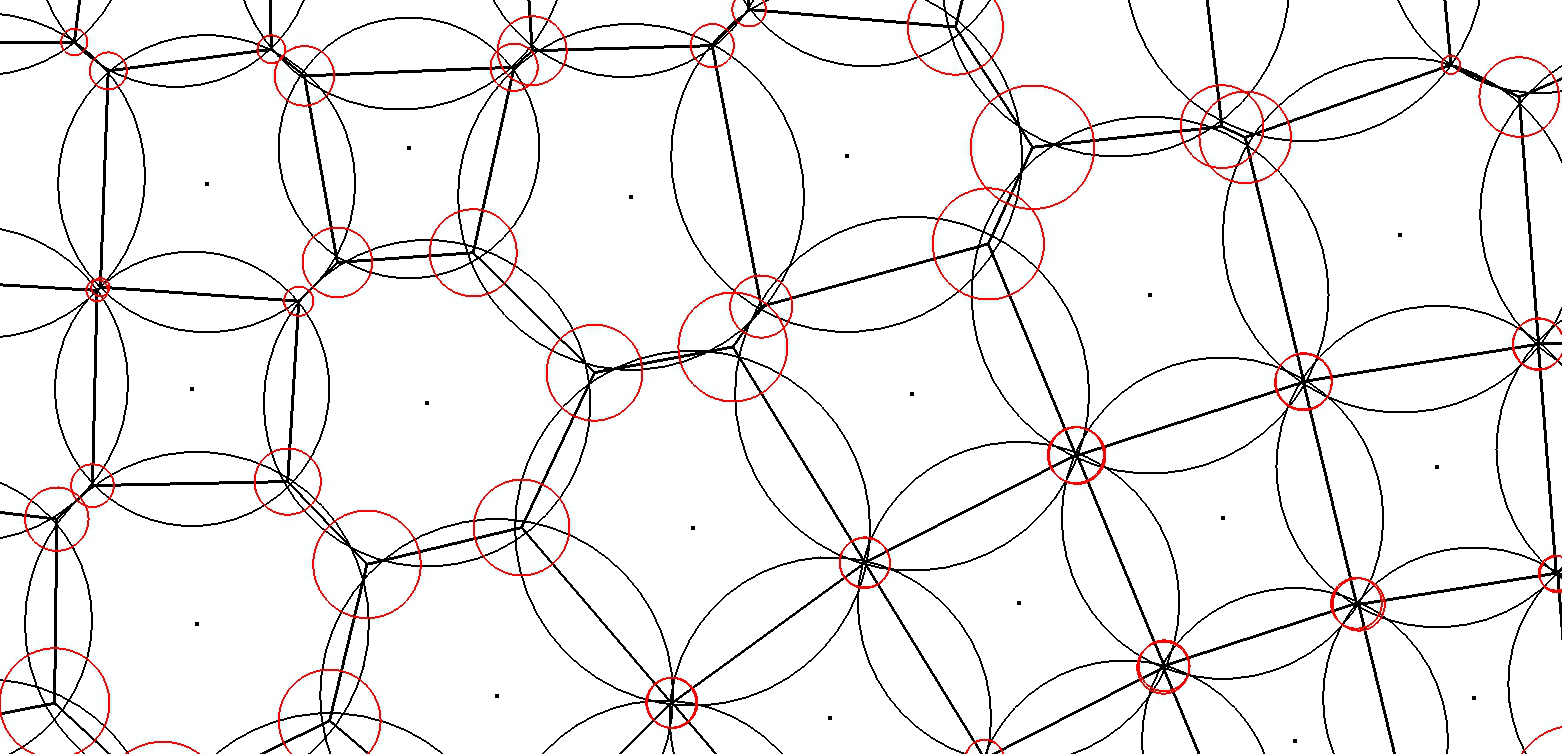}%
  }%
  \hspace{0.05\textwidth}%
  \subfloat[\label{fig:Delone-power:d}]{%
    \includegraphics[width=0.46\textwidth]{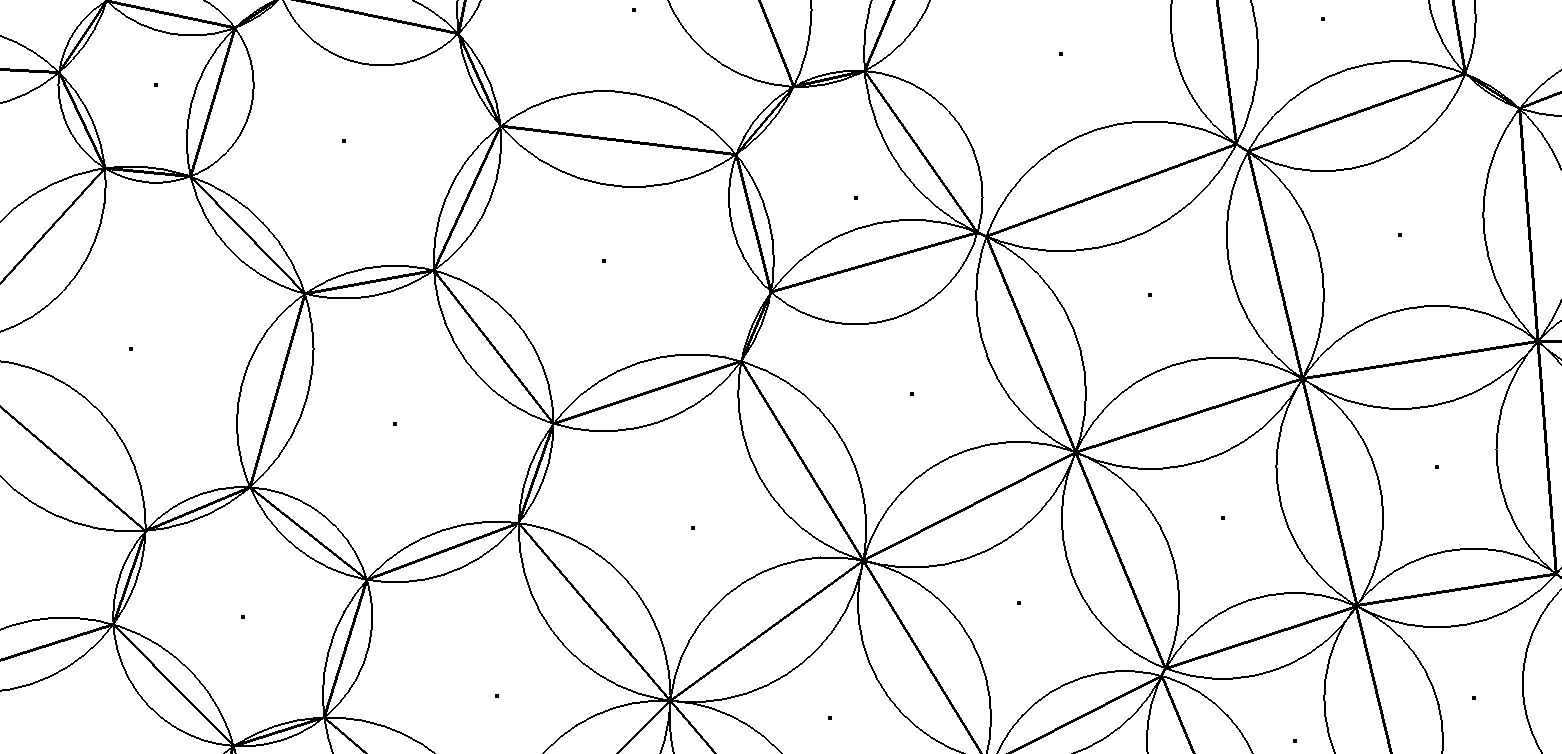}%
  }%
  \caption{Sequence of radical partitions
    \subref{fig:Delone-power:a}--\subref{fig:Delone-power:d};
   red circles correspond to the absolute values of powers}%
  \label{fig:Delone-power}%
\end{figure}

\Cref{fig:Delone-power} shows the evolution of the mesh fragment with added orthocircles, i.e., the circles centered at the dual vertices $\vb*{v}_k$ with $\sqrt{|\tau(\vb*{v}_k)|}$ as radii.
In the case $\tau(\vb*{v}_k) < 0$, the red circles are no longer real orthocircles; they are introduced just to evaluate the deviation of the radical cells from the Delaunay cells.

\begin{figure}[H]%
  \includegraphics[width=0.48\textwidth]{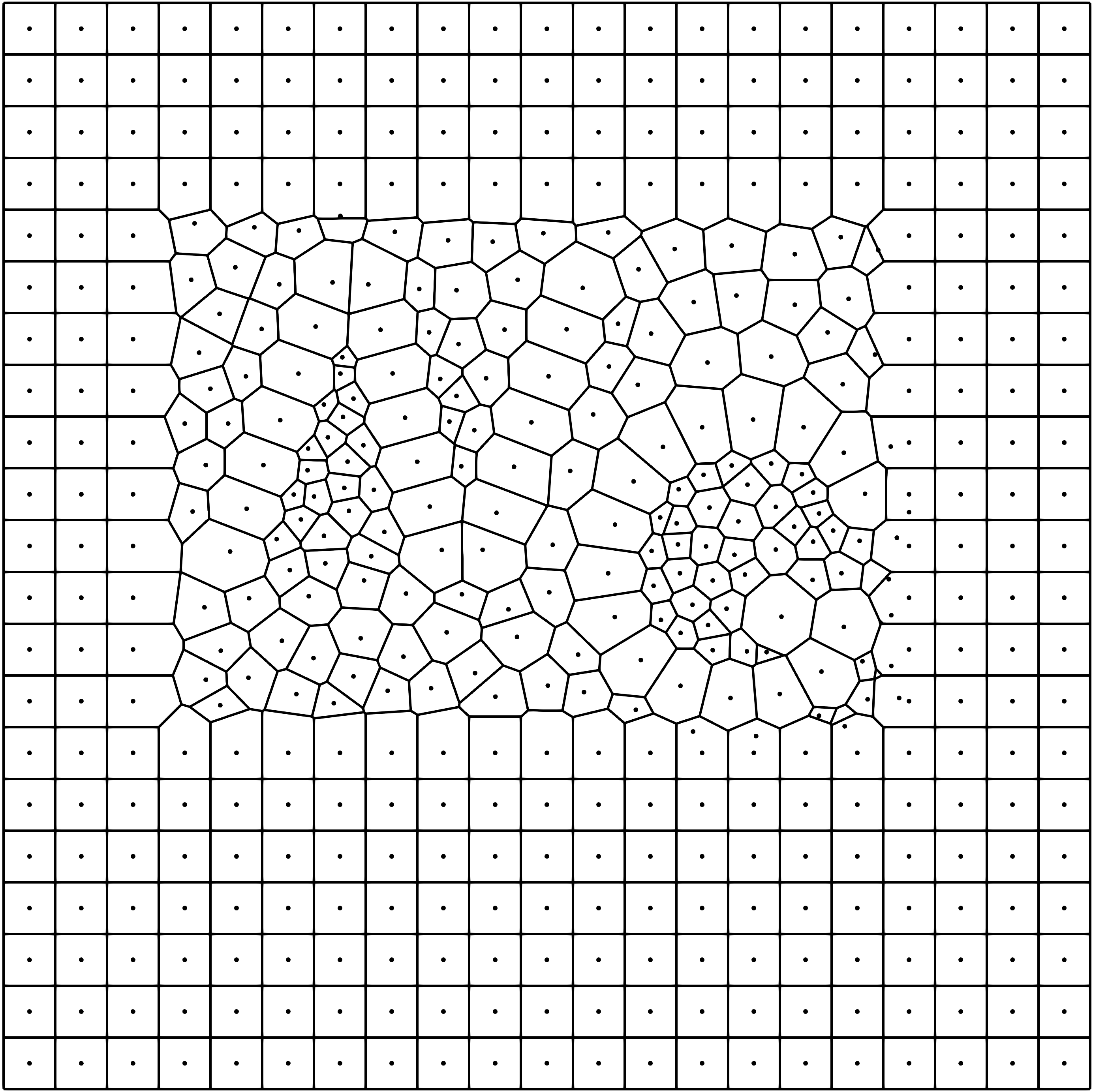}%
  \hfill{}%
  \includegraphics[width=0.48\textwidth]{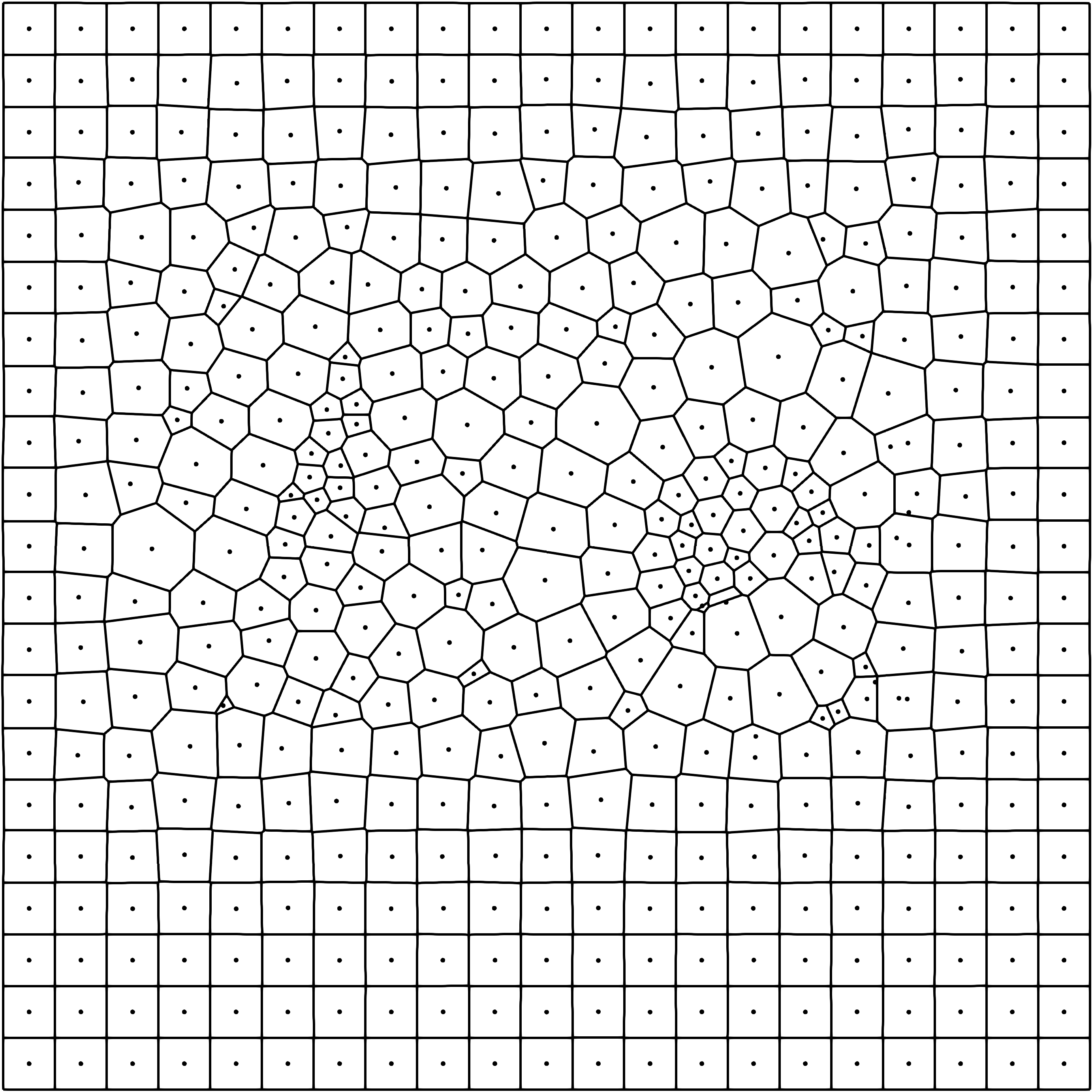}%
  \caption{NG logo: initial radical partition $\mathcal{R}_0$ and the final Delaunay partition $\mathcal{T}$}%
  \label{fig:ng:partition}%
\end{figure}

Next, we use the same algorithm to build the logo of the NUMGRID-2020 conference~\cite{Springer2021}.
The abbreviation NG is created using a set of fixed protecting circles incorporated into the square lattice of circles; a quasi-random set of circles is scattered around the letters.
The initial radical partition and the converged solution (Delaunay partition) are shown in \cref{fig:ng:partition}.

\begin{figure}[H]%
  \includegraphics[width=0.48\textwidth]{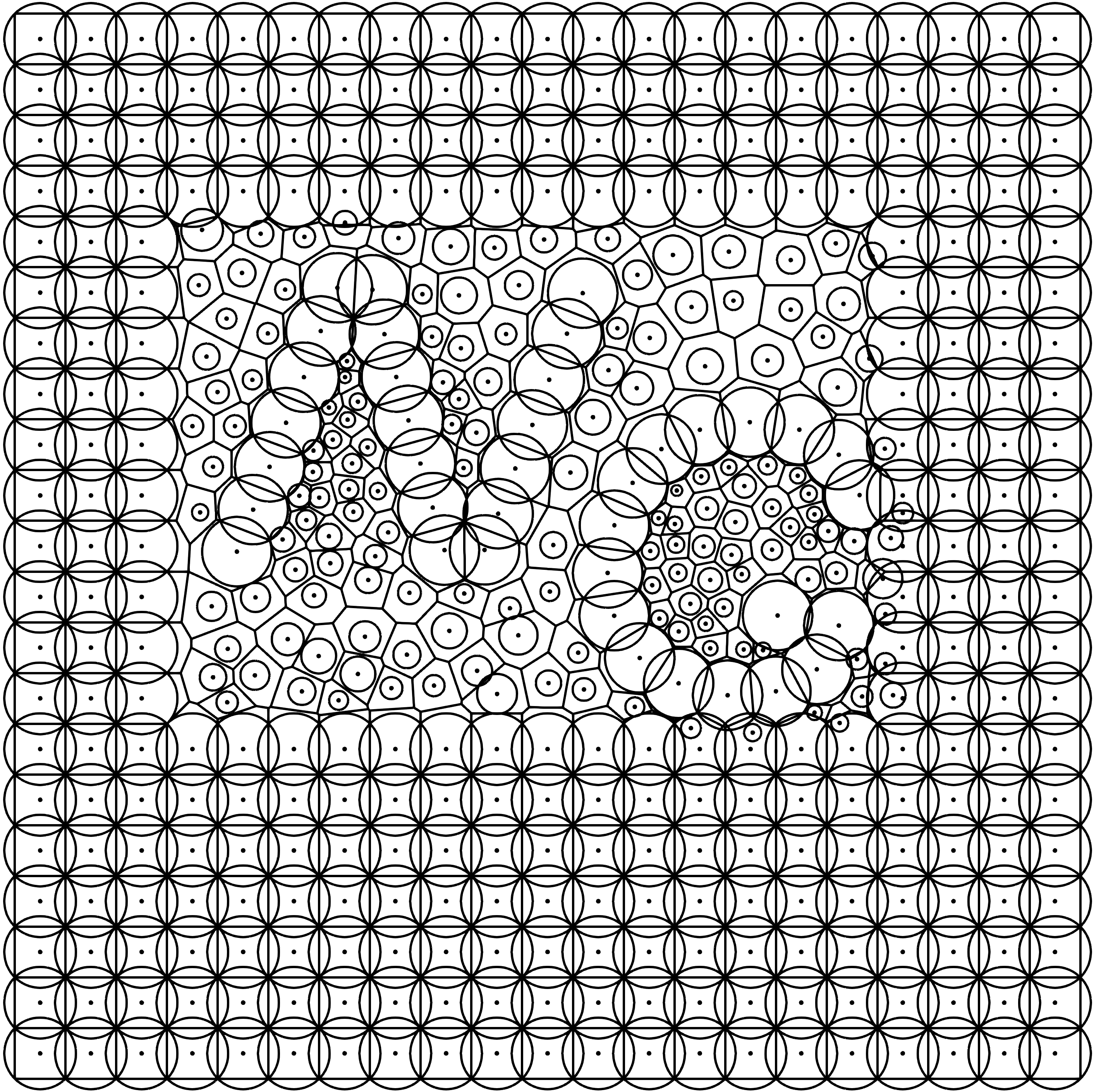}%
  \hfill{}%
  \includegraphics[width=0.48\textwidth]{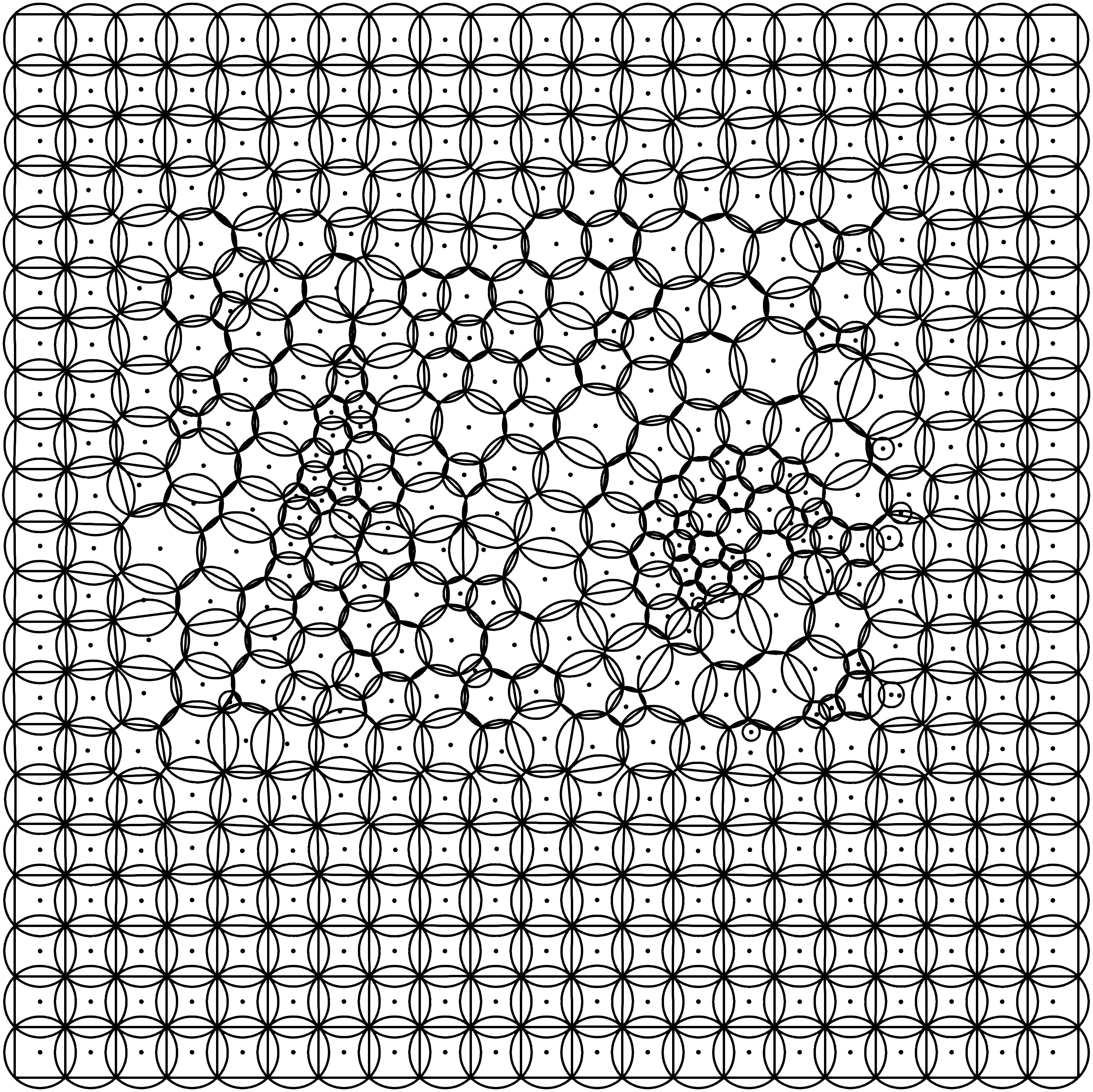}%
  \caption{NG logo: initial radical partition $\mathcal{R}_0$ with the initial set of circles and the final Delaunay partition $\mathcal{T}$ with Delaunay circles}%
  \label{fig:ng:partition:circle}%
\end{figure}

\Cref{fig:ng:partition:circle} shows radical partitions and the circles $B_i$ that generate them.
Note, that the final circles are actually Delaunay circles.

\begin{figure}[H]%
  \includegraphics[width=0.48\textwidth]{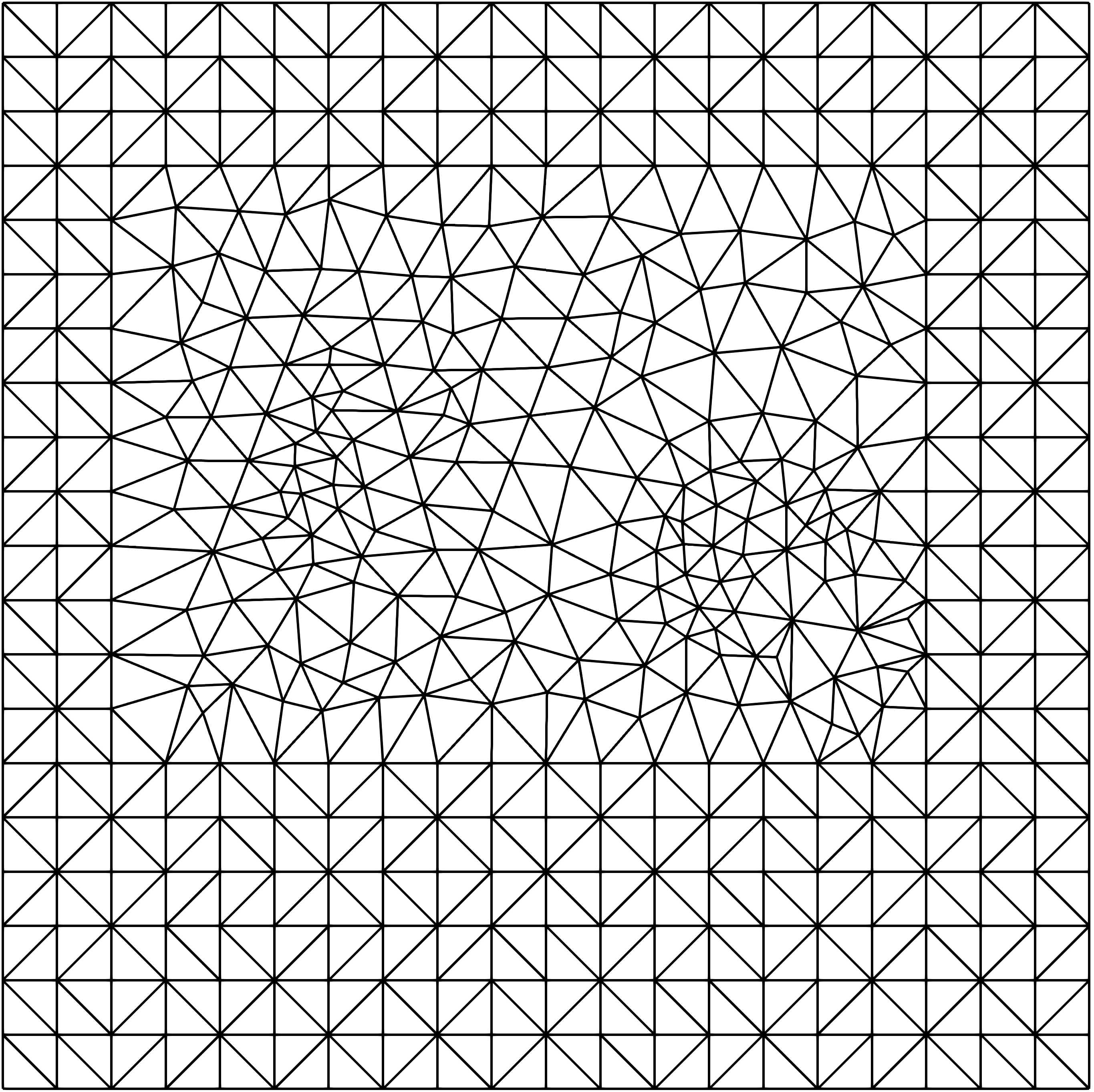}%
  \hfill{}%
  \includegraphics[width=0.48\textwidth]{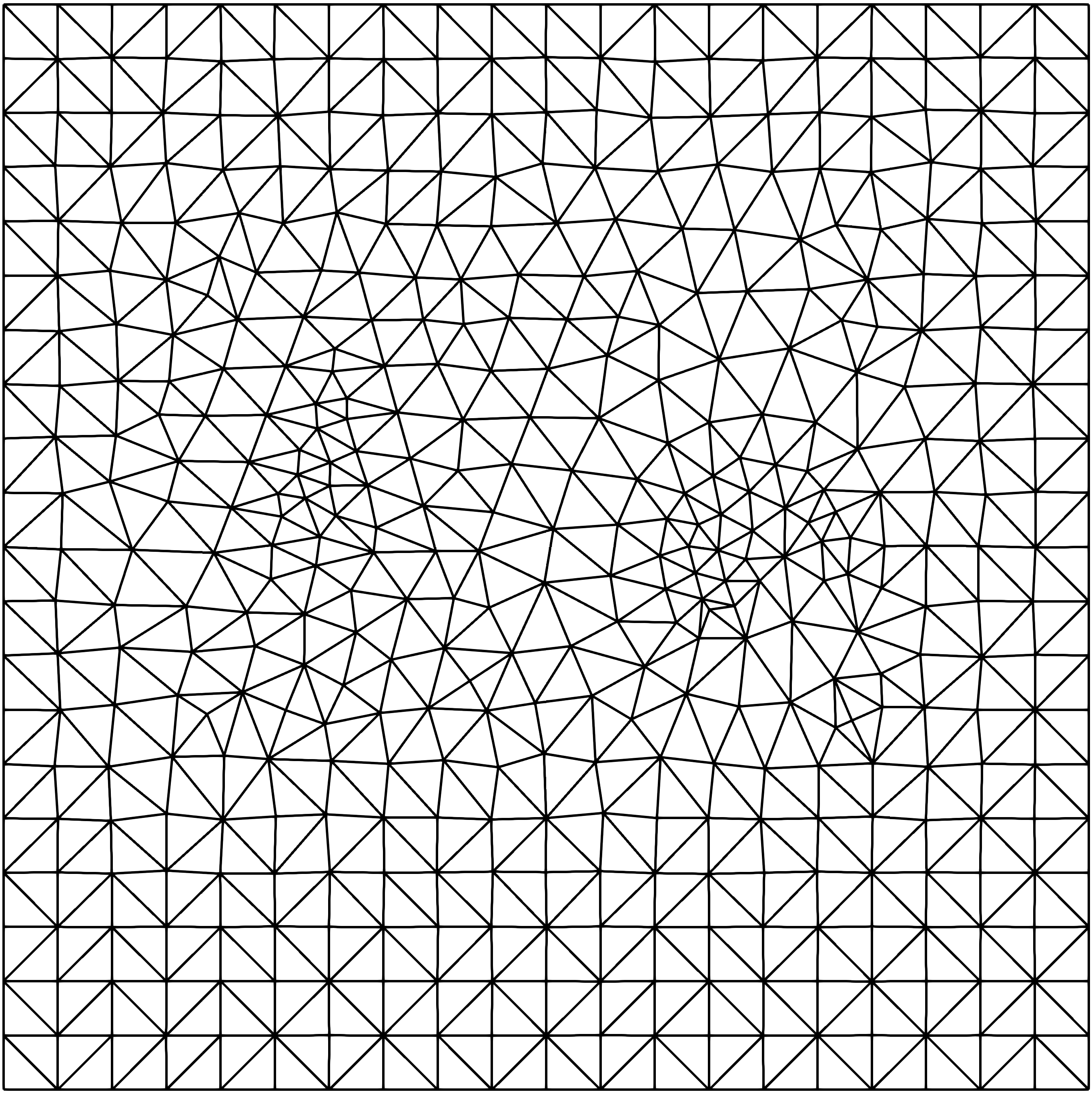}%
  \caption{NG logo: initial weighted Delaunay mesh $\mathcal{W}_0$ and the final Voronoi triangulation}%
  \label{fig:ng:mesh}%
\end{figure}

\Cref{fig:ng:mesh} shows the initial weighted Delaunay triangulation and the final Voronoi triangulation.

\begin{figure}[H]%
  \includegraphics[width=0.48\textwidth]{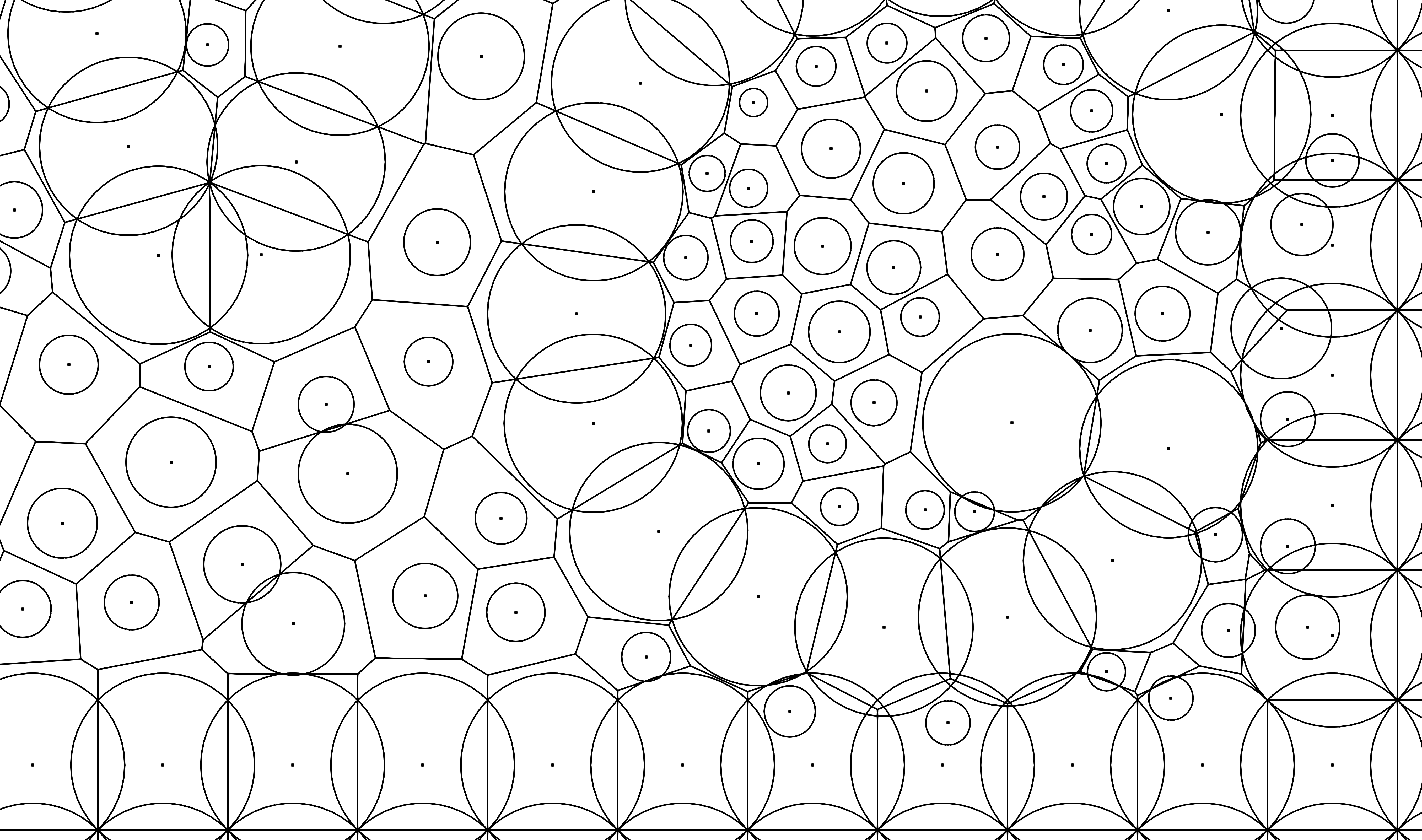}%
  \hfill{}%
  \includegraphics[width=0.48\textwidth]{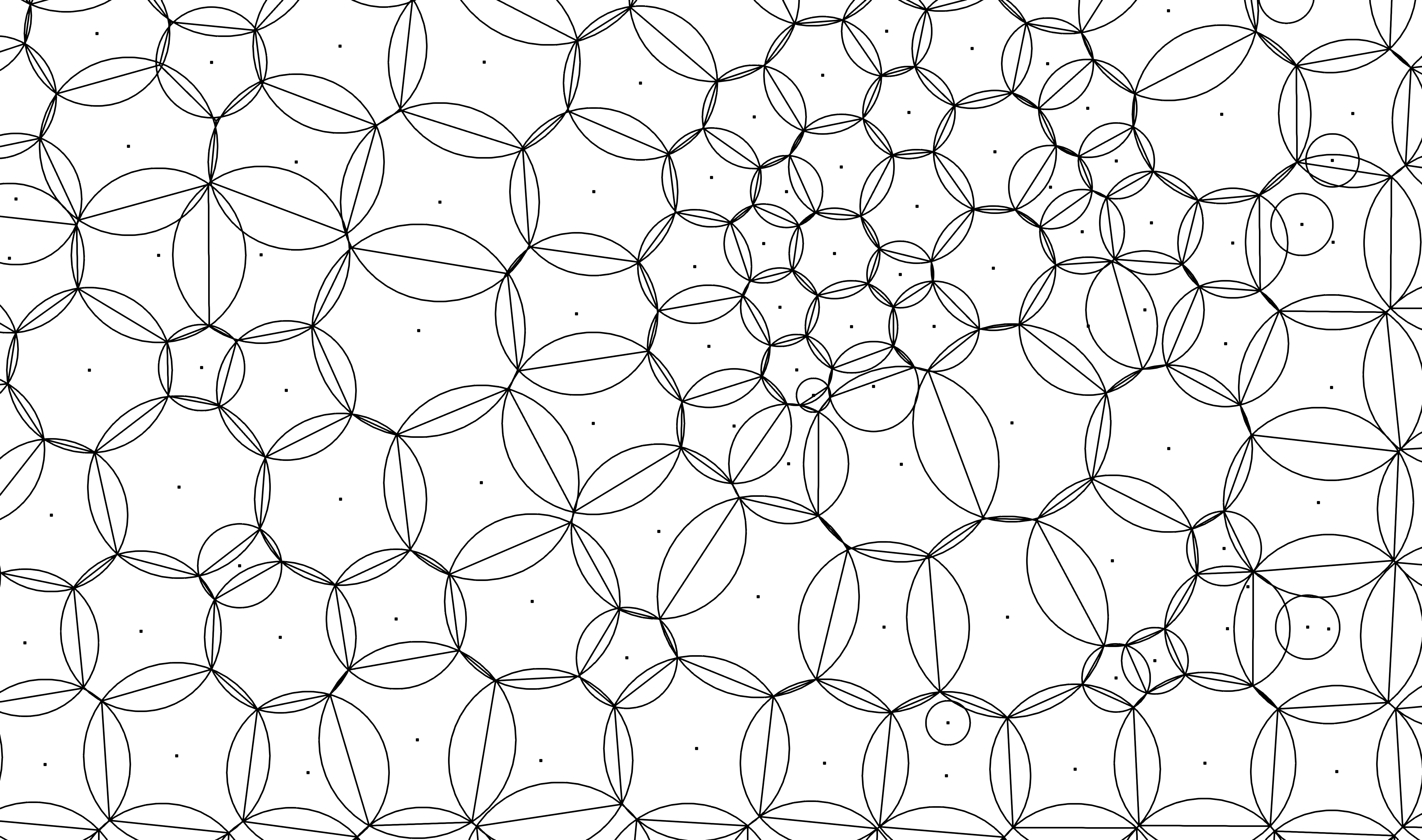}%
  \caption{NG logo: fragments of the initial radical and the final Delaunay partitions}%
  \label{fig:ng:partition:fragment}%
\end{figure}

Finally, \cref{fig:ng:partition:fragment} shows enlarged fragments of the initial radical partition and the converged radical partition, which coincides with the Delaunay partition.

\section{Discussion}
In the two-dimensional case, we have shown numerically that a radical partition can evolve into the polygonal Delaunay partition via evolution of the set of circles.
The problem setting is multi-dimensional, hence, it is expected that the algorithm is applicable in the three-dimensional case as well.
Although there is no existence result for this problem, we do not consider this as a crucial drawback.
As soon as the addition of new balls becomes an admissible operation, the existence result become trivial, because the projected function $\tilde{v}^*$ in each iteration is the solution of the problem.
However, the problem of the minimal addition of balls in order to construct a solution is an open one.
Numerical experiments suggest that the ability to locally add new circles/balls can be also important in order to attain a given mesh quality in the presence of constraints.

Note, that the concept of lifting points on the paraboloid allows for a development of nontrivial computational algorithms.
In~\cite{Alexa-2020}, the balls with a~priori unknown radii were assigned to the vertices of the surface triangulation in order to solve the recovery problem for the weighted Delaunay tetrahedralization of a point set matching prescribed boundary triangles.
In~\cite{Alexa-2020}, the problem of existence of such tetrahedralization is not solved.
Evidently, it is possible to define such a surface triangulation so that some of its faces cannot be weighted Delaunay faces.
In this case, the surface mesh should be refined: new vertices are added, which  alleviates the existence problem.

The idea of using the value of the discrete Dirichlet functional for an arbitrary lifting (for the roughness measure of the lifted surface) is not new.
In~\cite{Rippa-1990}, the minimal roughness principle for two-dimensional Delaunay triangulations was established.
At the moment, it is not clear how the minimal roughness property can be related to the presented results.

In practice, the mesh quality functionals should be optimized by using the manifold $\nabla F = 0$ as a constraint.
Obvious quality requirements are related to the elimination of small Delaunay and Voronoi edges and faces, which is subject of ongoing research.
In order to build good Delaunay-Voronoi meshes, one has to follow a given local sizing function, eliminate small Delaunay edges/faces and circles/balls, eliminate small Voronoi edges/faces, and give preference to the polyhedral rather than the simplicial Voronoi mesh.

\end{document}